\documentclass[11pt]{article}

\usepackage{amssymb}
\usepackage{amsmath}
\usepackage{theorem}
\usepackage{epsfig}
\usepackage{verbatim}
\usepackage{graphicx}
\usepackage{indentfirst}
\newtheorem{thm}{Theorem}[section]

\newtheorem{prop}[thm]{Proposition}

\newtheorem{rem}[thm]{Remark}

\newtheorem{lemma}[thm]{Lemma}

\newcommand{\R}{\Bbb{R}}

\newcommand{\T}{\mathbb{T}}
\newcommand{\D}{\displaystyle}

\newcommand{\sign}{{\rm sign}\thinspace}

\newcommand{\di}{{\rm div}\thinspace}
\newcommand{\curl}{{\rm curl}\thinspace}

\newcommand{\grad}{\nabla}

\newcommand{\dt}{\frac{d}{dt}}
\newcommand{\F}{\mathcal{F}}

\newcommand{\dxu}{\partial_{x_1}}
\newcommand{\dxd}{\partial_{x_2}}

\newcommand{\dpt}{\partial_t}
\newcommand{\dist}{{\rm dist}\thinspace}

\newcommand{\la}{\Lambda}

\newcommand{\al}{\alpha}
\newcommand{\ep}{\varepsilon}

\newcommand{\dpa}{\partial^{\bot}_{\alpha}}
\newcommand{\da}{\partial_{\alpha}}

\def\eop{{\ \vrule height 3pt width 3pt depth 0pt}}

\numberwithin{equation}{section}

\begin{document}

\author{Antonio C\'ordoba, Diego C\'ordoba and Francisco Gancedo}

\title{Interface evolution:  water waves  in 2-D.}
\date{Oct-29-08}

\maketitle

\begin{abstract}
We study the free boundary evolution between two irrotational, incompressible
and inviscid fluids in 2-D without surface tension.
 We prove local-existence in Sobolev
spaces when, initially, the difference of the gradients of the
pressure in the normal direction has the proper sign, an assumption
which is also known as the Rayleigh-Taylor condition. The well-posedness of  the full
 water wave problem was first obtained by Wu \cite{Wu}.
The methods introduced in this paper allows us to consider multiple cases:
with or without gravity, but also a closed boundary or
a periodic boundary with the fluids placed above and below it. It
is assumed that the initial interface does not touch itself, being a
part of the evolution problem to check that such property prevails
for a short time, as well as it does the Rayleigh-Taylor condition,
depending conveniently upon the initial data. The addition of the pressure equality to the contour dynamic
equations is obtained as a mathematical consequence, and not as a physical assumption, from the mere
fact that we are dealing with weak solutions of Euler's equation in the whole space.

\end{abstract}

\maketitle

\section{Introduction}

We consider the following evolution problem for the active scalar
$\rho = \rho(x,t)$, $x\in\mathbb{R}^2$, and $t\geq 0$:
\begin{equation}\label{conservacionmasa}
\rho_t + v\cdot\nabla\rho = 0,
\end{equation}
with a velocity $v = (v_1,v_2)$ satisfying the Euler equation
\begin{equation}\label{euler}
\rho(v_t+v\grad v)=-\nabla p-(0,\mathrm{g}\,\rho),
\end{equation}
and the incompressibility condition
\begin{equation}
\label{div0}\nabla\cdot v = 0.
\end{equation}
The free boundary is given by the discontinuity on the densities of
the fluids
\begin{equation*}\label{density}
\rho(x_1,x_2,t)=\left\{\begin{array}{cl}
                    \rho^1,& x\in\Omega^1(t)\\
                    \rho^2,& x\in\Omega^2(t)=\mathbb{R}^2 - \Omega^1(t),
                 \end{array}\right.
\end{equation*} where $\rho^1\neq\rho^2$ are constants.

We shall assume also that each fluid is irrotational, i.e. $\omega =
\nabla \times u = 0$, in the interior of each domain $\Omega^j$
($j=1,2$). The main purpose of this paper is to understand the
evolution of the free boundary, but we shall also take the point of
view of having weak solutions in the whole space presenting a
discontinuity in the density along the interface. Under the
hypothesis that at the initial time we have smooth velocity fields
$v^1$, $v^2$ whose values at the interface differs only in the
tangential direction it follows that, for a certain time $t>0$, the
vorticity  $\omega$ will be supported on the free boundary curve
$z(\al,t)$ and it has the form
\begin{equation*}
\omega(x,t)=\varpi(\al,t)\delta(x-z(\al,t)).
\end{equation*}
Here we shall consider two types of geometries, namely periodicity
in the horizontal space variable, says  $z(\al+ 2k\pi,t) =
z(\al,t)+(2k\pi,0)$, or the case of a closed contour $z(\al+
2k\pi,t) = z(\al,t)$. We shall assume also that we have infinite
depth. In \cite{Lannes} fluids of finite depth were considered.

In section 2 our first step will be to show the equality of pressure
at each side of the free boundary, when we understand the system
(\ref{conservacionmasa}--\ref{div0})  in a weak sense (see
Proposition \ref{presiones}).

The free boundary $z(\alpha,t)$ evolves with a velocity field coming from
Biot-Savart law, which can be explicitly computed and it is given by
the Birkhoff-Rott integral of the amplitude $\varpi$ along the
interface curve:
\begin{align}
\begin{split}\label{fibr}
BR(z,\varpi)(\al,t)=\frac{1}{2\pi}
PV\int \frac{(z(\al,t)-z(\beta,t))^{\bot}}{|z(\al,t)-
z(\beta,t)|^2}\varpi(\beta,t)d\beta,
\end{split}
\end{align}
where $PV$ denotes principal value \cite{St3}.
It gives us the velocity field at the interface to which we can
subtract any term in the tangential direction without modifying the
geometric evolution of the curve
\begin{align}
\begin{split}\label{fullpm}
z_t(\al,t)&=BR(z,\varpi)(\al,t)+c(\al,t)\da z(\al,t).
\end{split}
\end{align}
A wise choice of $c(\al,t)$ namely:
\begin{align}
\begin{split}\label{fc}
c(\al,t)&=\frac{\al+\pi}{2\pi}\int_{-\pi}^\pi \frac{\da z(\al,t)}{|\da
z(\al,t)|^2}\cdot \da BR(z,\varpi)(\al,t) d\al\\
&\quad-\int_{-\pi}^\al \frac{\da z(\beta,t)}{|\da
z(\beta,t)|^2}\cdot \partial_{\beta} BR(z,\varpi)(\beta,t)d\beta,
\end{split}
\end{align}
allows us to accomplish the fact that the length of the tangent
vector to $z(\al,t)$ be just a function in the variable $t$ only \cite{Hou}:
$$
A(t)=|\da z(\al,t)|^2.
$$

Then we can close the system using Bernoulli's law with the
equation:
\begin{align}
\begin{split}\label{fw}
\varpi_t(\al,t)&=-2A_\rho \partial_t BR(z,\varpi)(\al,t)\cdot
\da z(\al,t)-A_\rho\da( \frac{|\varpi|^2}{4|\da z|^2})(\al,t) +\da (c\, \varpi)(\al,t)\\
&\quad +2A_{\rho}c(\al,t)\da BR(z,\varpi)(\al,t)\cdot\da
z(\al,t)-2A_{\rho}\mathrm{g}\da z_2(\al,t),
\end{split}
\end{align}
where $$A_\rho=\frac{\rho_2-\rho_1}{\rho_2+\rho_1}$$ is the Atwood
number.

We shall use the notation $T$ for the following operator (depending
on the curve $z(\al,t)$) acting on $u(\al,t)$ by the formula
\begin{equation}\label{dti}
T(u)(\al,t)=2BR(z,u)(\al,t)\cdot\da z(\al,t).
\end{equation}
The inversibility of $(I+A_\rho T)$ (see \cite{BMO}) allows us to
write the equation  \eqref{fw} in the following more convenient explicit
manner:
\begin{equation}\label{iopti}
\varpi_t(\al,t)=(I+A_\rho T)^{-1}(A_\rho R(z,\varpi)+\da(c\varpi))(\al,t).
\end{equation}

Next let us give the function which measures the arc-chord condition \cite{Y}
\begin{equation}\label{df}
\mathcal{F}(z)(\alpha,\beta,t)=\frac{|\beta|}{|z(\alpha,t)-z(\alpha-\beta,t)|}\qquad
\forall\,\alpha,\beta\in(-\pi,\pi),
\end{equation}
and $$\mathcal{F}(z)(\alpha,0,t)=\frac{1}{|\da z(\alpha,t)|}.$$

Finally  following references \cite{BHL} and
\cite{AM} we introduce the auxiliary function $\varphi(\al,t)$ which will allow us  to integrate the evolution equation
\begin{equation}\label{fphi}
\varphi(\al,t)=\frac{\varpi(\al,t)}{2|\da z(\al,t)|}-c(\al,t)|\da
z(\al,t)|.
\end{equation}

Our main result consists on local existence for the water wave
problem: $\rho_1=0$. We prove that there is  a positive time
$\mathrm{T}$ (depending upon the initial condition) for which there
exists  a solution of the equations (\ref{fibr}--\ref{fw}) with
$\rho_1=0$ during the  time interval $[0,\mathrm{T}]$ so long as the
initial data satisfy $z_0(\alpha)\in H^{k}$, $\varphi_0(\al)\in
H^{k-\frac12}$ and $\varpi_0(\al)\in H^{k-1}$ for $k\geq 4$,
$\F(z_0)(\al,\beta)<\infty$,  and
$$\sigma_0(\al)= -(\nabla p^2(z_0(\alpha),0) - \nabla
p^1(z_0(\alpha),0))\cdot\partial^{\perp}_{\alpha}z_0(\alpha) > 0,$$
where $p^j$ denote the pressure in $\Omega^j$.

\begin{thm}\label{theorem2}
Let $z_0(\alpha)\in H^{k}$, $\varphi_0(\alpha)\in
H^{k-\frac12}$ and $\varpi_0(\al)\in H^{k-1}$ for $k\geq 4$, $\F(z_0)(\al,\beta)<\infty$, and
$$\sigma_0(\al)=-(\nabla p^2(z_0(\alpha),0) - \nabla
p^1(z_0(\alpha),0))\cdot\partial^{\perp}_{\alpha}z_0(\alpha)
> 0.$$
Then there exists a time $\mathrm{T}>0$ so that we have a solution to
(\ref{fibr}--\ref{fw}) in the case  $\rho_1=0$, where $z(\al,t)\in
C^{1}([0,\mathrm{T}];H^k)$ and $\varpi(\al,t)\in
C^{1}([0,\mathrm{T}];H^{k-1})$ with $z(\al,0)=z_0(\al)$ and
$\varpi(\al,0)=\varpi_0(\alpha)$.
\end{thm}

The first results concerning the Cauchy problem for the linearized
version in Sobolev spaces are due to \cite{Craig}, \cite{Nalinov}
and \cite{Yoshira}. In her important work \cite{Wu} (see also
\cite{Wu2}) S. Wu was able to prove that the presence of the
gravitational field,
 together with the hypothesis about the asymptotic flatness of the fluid domains, implies that
 %if the Rayleigh-Taylor signum
%condition is satisfied at the initial time then it must continue to hold at any
%later time, so long as the interface is well-defined
the Rayleigh-Taylor signum condition must hold so long as the
interface is well-defined. In our treatment we can also get local
solvability even in the absence of gravity, or for a closed
contour,
 whenever the Rayleigh-Taylor and the arc-chord conditions are initially satisfied.

Besides the significant work of S. Wu that has been referred before,
we can also quote the interesting paper \cite{AM} where they get
energy estimates on the free boundary and the amplitude of the
vorticity,
 under the time dependent assumption of the arc-chord property. These authors make also use of the fact
 obtained by Wu about the persistence of the Rayleigh-Taylor sign condition.

In our approach the explicit control upon the evolution of the arc-chord relation of the free
boundary is especially emphasized, together with the inversion of the operator $(I+T)$,
which gives us the equation for the time derivative of the vorticity amplitude in terms of
 the curve (see equations (\ref{dti}--\ref{iopti}) with $\rho_1=0$). The architecture of our
  proof relies upon different energy estimates for the quantities involved (Sobolev norms
   for $z$, $\varpi$, arc-chord and Rayleigh-Taylor condition). But in order to fix together
    its different parts it becomes crucial to get explicit upper bounds on the operator
    $(I+T)^{-1}$ on different Sobolev spaces. Here we continue the method introduced in
    \cite{DY2} and \cite{DY3}, where conformal mappings, Hopf maximum principle and Dahlbert-Harnack
    inequality up to the boundary, for nonnegative harmonic functions, play a central role.

 In the following interesting works by Christodoulou-Lindblad \cite{Chris},
  Lindblad \cite{Lindblad}, Coutand-Shkoller \cite{Shkoller}, Shatah-Zeng \cite{Shatah}
   and Zhang-Zhang \cite{ZZ} the rotational case have been also considered.
   Let us point out that the evolution of the sign of Rayleigh-Taylor condition
   is crucial in our proof \cite{DY3}, because it allows to get rid of the highest
   order derivatives in the evolution equation of the Sobolev norms of the curve (section 8).

%%%%%%%%%%%%%%%%%%%%%%%%%%%%%%%%%%%%%%%%%%%%%%%%%%%%%%%%%%%%%%%%%%%%%%%%%%%%%%%%%%%%%%%%%%%%%%%%%%%%%%%%%%%%%%%%%%%%%%%%%%
%%%%%%%%%%%%%%%%%%%%%%%%%%%%%%%%%%%%%%%%%%%%%%%%%%%%%%%%%%%%%%%%%%%%%%%%%%%%%%%%%%%%%%%%%%%%%%%%%%%%%%%%%%%%%%%%%%%%%%%%%%

\section{The evolution equation}

%%%%%%%%%%%%%%%%%%%%%%%%%%%%%%%%%%%%%%%%%%%%%%%%%%%%%%%%%%%%%%%%%%%%%%%%%%%%%%%%%%%%%%%%%%%%%%%%%%%%%%%%%%%%%%%%%%%%%%%%%%
%%%%%%%%%%%%%%%%%%%%%%%%%%%%%%%%%%%%%%%%%%%%%%%%%%%%%%%%%%%%%%%%%%%%%%%%%%%%%%%%%%%%%%%%%%%%%%%%%%%%%%%%%%%%%%%%%%%%%%%%%%

We shall consider weak solutions of the system
(\ref{conservacionmasa}--\ref{div0}); that is for any smooth
functions $\zeta$, $\eta$ and $\chi$, compactly supported on
$[0,\mathrm{T})\times\R^2$ i.e. lying in the space
$C^\infty_c([0,\mathrm{T})\times\R^2)$, we have
\begin{equation}\label{weakconservacionmasa}
\int_0^\mathrm{T}\int_{\R^2}\rho(\zeta_t + v\cdot\nabla\zeta)
dxdt+\int_{\R^2}\rho_0(x)\zeta(x,0)dx = 0,
\end{equation}
\begin{equation}\label{weakeuler}
\int_0^\mathrm{T}\int_{\R^2}\big(\rho v\cdot (\eta_t + v\cdot\nabla\eta)
+p\nabla\cdot\eta-(0,\mathrm{g}\rho)\cdot\eta\big)dxdt+\int_{\R^2}\rho_0(x)v_0(x)\cdot\eta(x,0)dx = 0,
\end{equation}
and
\begin{equation}\label{weakdiv0}
\int_0^\mathrm{T}\int_{\R^2} v\cdot \nabla \chi dxdt= 0.
\end{equation}

Here $\rho$ is defined by
\begin{equation}\label{rhocomoes}
\rho(x_1,x_2,t)=\left\{\begin{array}{cl}
                    \rho^1,& x\in\Omega^1(t)\\
                    \rho^2,& x\in\Omega^2(t),
                 \end{array}\right.
\end{equation} where $\rho^1\neq\rho^2$. It is assumed that the vorticity is given by a
delta function on the curve $\partial\Omega^j(t)$ multiplied by an
amplitude and has the form
\begin{equation}\label{omegacomoes}
\omega(x,t)=\varpi(\al,t)\delta(x-z(\al,t)).
\end{equation}

 Then using the Biot-Savart law we get
\begin{equation}\label{vfc}
v(x,t)=\frac{1}{2\pi}PV\int\frac{(x-z(\beta,t))^\bot}{|x-z(\beta,t)|^2}
\varpi(\beta,t)d\beta
\end{equation}
for $x$ not lying on the curve $z(\al,t)$, and
\begin{align}\label{limitevelocidad}
\begin{split}
v^2(z(\al,t),t)&=BR(z,\varpi)(\al,t)+\frac12\frac{\varpi(\al,t)}{|\da
z(\al,t)|^2}\da z(\al,t),\\
v^1(z(\al,t),t)&=BR(z,\varpi)(\al,t)-\frac12\frac{\varpi(\al,t)}{|\da
z(\al,t)|^2}\da z(\al,t),
\end{split}
\end{align}
where $v^j(z(\alpha,t),t)$ denotes the limit velocity field obtained
approaching the boundary in the normal direction inside $\Omega^j$
and $BR(z,\varpi)(\al,t)$ is given by \eqref{fibr}. It is easy to
check that \eqref{weakdiv0} is satisfied by $v$ given as in
\eqref{vfc}. Furthermore, we have that the identity of the weak
formulation \eqref{weakconservacionmasa} is verified so long as the
following equality holds (see \cite{DY}):
\begin{equation}\label{vnormal}
z_t(\al,t)\cdot \dpa z(\al,t)=BR(z,\varpi)(\al,t)\cdot\dpa z(\al,t).
\end{equation}

\begin{prop}\label{presiones}
Let us consider a weak solution $(\rho,v,p)$ satisfying (\ref{weakconservacionmasa}--\ref{weakdiv0})
where $\rho$ is given by $\eqref{rhocomoes}$ and $\curl v=\omega$ by \eqref{omegacomoes}.
 Then we have the following identity $$p^1(z(\al,t),t)=p^2(z(\al,t),t),$$
where $p^j(z(\alpha,t),t)$ denotes the limit pressure obtained
approaching the boundary in the normal direction inside $\Omega^j$.
\end{prop}

Proof: We shall show that the Laplacian of the pressure is as follows
$$
\Delta p(x,t)= F(x,t)+f(\alpha,t)\delta(x-z(\al,t)),
$$
where $F$ is regular in $\Omega^j(t)$ although discontinuous on
$z(\al,t)$, and the amplitude of the Dirac distribution $f$ is
regular. Then the inverse of the Laplacian by means of the
Newtonian potential gives the continuity of the pressure on the
free boundary (see \cite{DY2}).

We also shall use an ad hoc integration by parts for the
derivatives of the velocity. The expression for the conjugate of
the velocity in complex variables
$$
\overline{v}(z,t)=\frac{1}{2\pi i}PV\int \frac{1}{z-z(\al,t)}\varpi(\al,t)d\al,
$$
for $z\neq z(\al,t)$ allows us to accomplish the fact that
\begin{align*}
\partial_z\overline{v}(z,t)&=\frac{1}{2\pi i}PV\int \frac{-\varpi(\al,t)}{(z-z(\al,t))^2}d\al=\frac{1}{2\pi i}
PV\int \frac{-\da z(\al,t)}{(z-z(\al,t))^2}\frac{\varpi(\al,t)}{\da z(\al,t)}d\al,
\end{align*}
and therefore
\begin{equation}\label{dzv}
\partial_z\overline{v}(z,t)=\frac{1}{2\pi i}PV\int \frac{1}{z-z(\al,t)}\da(\frac{\varpi}{\da z})(\al,t) d\al,
\end{equation}
for a regular parametrization with $\da z(\al,t)\neq 0$. In a similar way
\begin{equation}\label{vt}
\overline{v}_t(z,t)=\frac{1}{2\pi i}PV\int \frac{1}{z-z(\al,t)}\varpi_t(\al,t)d\al-\frac{1}{2\pi i}
PV\int \frac{1}{z-z(\al,t)}\da(\frac{z_t\varpi}{\da z})(\al,t) d\al,
\end{equation}
and
\begin{align}\label{d2zv}
\partial^2_z\overline{v}(z,t)&=\frac{1}{2\pi i}PV\int \frac{1}{z-z(\al,t)}\da(\frac{1}{\da z}\da(\frac{\varpi}{\da z}))(\al,t) d\al.
\end{align}
These identities help us to get the values of $\grad
v^j(z(\al,t),t)$, $v^j_t(z(\al,t),t)$ and $\grad^2 v^j(z(\al,t),t)$
which are obtained as limits approaching the boundary in the normal
direction inside $\Omega^j(t)$.

To get the stated formula for the pressure we start with identity
 \eqref{weakeuler} choosing $\eta(x,t)=\grad\lambda(x,t)$. Then
\begin{align*}
\int_0^\mathrm{T}\!\!\int_{\R^2}p\Delta\lambda dxdt&=\int_0^\mathrm{T}\!\!\int_{\R^2}(0,\mathrm{g}\rho)\cdot\grad\lambda dxdt-
\int_0^\mathrm{T}\!\!\int_{\R^2}\rho v\cdot \grad\lambda_t dxdt\\
&\quad-\int_0^\mathrm{T}\!\!\int_{\R^2} \rho v\cdot(v\cdot\grad^2\lambda)dxdt-\int_{\R^2}\rho_0(x)v_0(x)\cdot\grad\lambda(x,0)dx\\
&=I_1+I_2+I_3+I_4.
\end{align*}
Let us define
$\Omega_{\ep}^1(t)=\{x\in\Omega^1(t):\dist(x,\partial\Omega^1(t))\geq\ep\}$
and
$\Omega_{\ep}^2(t)=\{x\in\Omega^2(t):\dist(x,\partial\Omega^2(t))\geq\ep\}$,
we have
\begin{align*}
I_1&=\lim_{\ep\rightarrow 0}
\int_0^\mathrm{T}\!\!\int_{\Omega^1_{\ep}(t)}\mathrm{g}\rho^1\dxd\lambda dxdt+
\int_0^\mathrm{T}\!\!\int_{\Omega^2_{\ep}(t)}\mathrm{g}\rho^2\dxd\lambda dxdt\\
&=\int_0^\mathrm{T}\!\!\int_{-\pi}^\pi(\rho^2-\rho^1)\mathrm{g}\da
z_1(\al,t) \lambda(z(\al,t),t)d\al dt,
\end{align*}
and we can consider the term $(\rho^2-\rho^1)\mathrm{g}\da
z_1(\al,t)$ as being part of the function $f(\al,t)$.

Regarding the term $I_2$ we integrate by parts in the variable $t$
to obtain
\begin{align*}
I_2&=\lim_{\ep\rightarrow 0}
-\int_0^\mathrm{T}\!\!\int_{\Omega^1_{\ep}(t)}\rho^1v^1\cdot\grad\lambda_t dxdt-
\int_0^\mathrm{T}\!\!\int_{\Omega^2_{\ep}(t)}\rho^2v^2\cdot\grad\lambda_t dxdt\\
&=J_1+J_2-I_4
\end{align*}
where
$$
J_1=\int_0^\mathrm{T}\!\!\int_{\R^2}\rho v_t\cdot \grad\lambda dxdt,
$$
and
$$
J_2=\int_0^\mathrm{T}\!\!\int_{-\pi}^\pi(\rho^2v^2(z(\al,t),t)-\rho^1v^1(z(\al,t),t))\cdot
\grad\lambda(z(\al,t),t)\,z_t(\al,t)\cdot\dpa z(\al,t)d\al dt.
$$
In $J_1$ we use formula \eqref{vt} to get the limit on the boundary
of $v_t(x,t)$. Again we first integrate by parts in $J_1$ and then
take  the limit when $\epsilon\rightarrow 0$. Since in each
$\Omega_{\ep}^j(t)$ $v_t$ is regular and $\di v_t=0$, it follows
that
$$
J_1=\int_0^\mathrm{T}\!\!\int_{-\pi}^\pi (\rho^2 v^2_t(z(\al,t),t)-\rho^1
v^1_t(z(\al,t),t))\cdot\dpa z(\al,t)  \lambda(z(\al,t),t) dxdt.
$$
As before we may consider $(\rho^2 v^2_t(z(\al,t),t)-\rho^1
v^1_t(z(\al,t),t))\cdot\dpa z(\al,t)$ as being a part of $f(\al,t)$.

Next \eqref{limitevelocidad} yields the splitting  $J_2=K_1+K_2$ where
$$
K_1=\int_0^\mathrm{T}\!\!\int_{-\pi}^\pi(\rho^2-\rho^1)BR(z,\varpi)(\al,t)\cdot
\grad\lambda(z(\al,t),t)\,z_t(\al,t)\cdot\dpa z(\al,t)d\al dt,
$$
$$
K_2=\int_0^\mathrm{T}\!\!\int_{-\pi}^\pi(\rho^2+\rho^1)\frac{\varpi(\al,t)}{2|\da z(\al,t)|^2}\da z(\al,t)\cdot
\grad\lambda(z(\al,t),t)\,z_t(\al,t)\cdot\dpa z(\al,t)d\al dt.
$$
Integrating by parts in $\al$ we can write
$$
K_2=-\int_0^\mathrm{T}\!\!\int_{-\pi}^\pi(\rho^2+\rho^1)
\lambda(z(\al,t),t)\da(\frac{\varpi}{2|\da z|^2}z_t\cdot\dpa
z)(\al,t)d\al dt,
$$
giving us another term of $f(\al,t)$.

Let us introduce now  the decomposition  $I_3=J_3+J_4+J_5+J_6$ where
$$
J_3=-\int_0^\mathrm{T}\!\!\int_{\R^2}\rho(v_1)^2\dxu^2\lambda dxdt,\qquad
J_4=-\int_0^\mathrm{T}\!\!\int_{\R^2}\rho v_1v_2\dxd\dxu\lambda dxdt,
$$
$$
J_5=-\int_0^\mathrm{T}\!\!\int_{\R^2}\rho v_1v_2\dxu\dxd\lambda dxdt,\qquad
J_6=-\int_0^\mathrm{T}\!\!\int_{\R^2}\rho(v_2)^2\dxd^2\lambda dxdt.
$$
Using the sets $\Omega^j_\ep(t)$ and the identity \eqref{dzv} we get
\begin{align*}
J_3&=\int_0^\mathrm{T}\!\!\int_{\R^2}2\rho v_1\dxu v_1\dxu\lambda dxdt\\
&\quad+\int_0^\mathrm{T}\!\!\int_{-\pi}^\pi\big(\rho^2(v^2_1(z(\al,t),t))^2-\rho^1(v_1^1(z(\al,t),t))^2\big)\dxu\lambda(z(\al,t),t)\da z_2(\al,t) d\al dt\\
&=K_3+K_4.
\end{align*}
The term $K_3$ trivializes because  the ad hoc integration by parts
formula together with the identity \eqref{d2zv} gives
\begin{align*}
K_3&=-\int_0^\mathrm{T}\!\!\int_{\R^2}2\rho (v_1\dxu^2 v_1+(\dxu v_1)^2)\lambda dxdt
-\int_0^\mathrm{T}\!\!\int_{-\pi}^\pi \widetilde{f}(\al,t)\lambda(z(\al,t),t) d\al dt.
\end{align*}
where $\widetilde{f}(\al,t)=2(\rho^2v^2_1(z(\al,t),t)\dxu
v^2_1(z(\al,t),t)-\rho^1 v_1^1(z(\al,t),t) \dxu
v^1_1(z(\al,t),t))\da z_2(\al,t)$, and the first term in $K_3$ is
part of $F(x,t)$ while the second lies in $f(\al,t)$.

We can rewrite $K_4$ as follows
\begin{align}
\begin{split}\label{k4}
K_4&=(\rho^2\!-\!\rho^1)\int_0^\mathrm{T}\!\!\int_{-\pi}^\pi [(BR_1)^2+\frac{\varpi^2}{4}\frac{(\da z_1)^2}{|\da z|^4}]\dxu\lambda(z)\da z_2 d\al dt\\
&\quad+(\rho^2\!+\!\rho^1)\int_0^\mathrm{T}\!\!\int_{-\pi}^\pi \varpi BR_1\frac{\da z_1}{|\da z|^2}\dxu\lambda(z)\da z_2 d\al dt.
\end{split}
\end{align}
Next we continue analogously  with $J_4$
\begin{align*}
J_4&=\int_0^\mathrm{T}\!\!\int_{\R^2}\rho (v_2\dxd v_1 +v_1\dxd v_2)\dxu\lambda dxdt\\
&\quad-\int_0^\mathrm{T}\!\!\int_{-\pi}^\pi\big(\rho^2(v^2_1v^2_2)(z(\al,t),t)-\rho^1(v_1^1v_2^1)(z(\al,t),t)\big)\dxu\lambda(z(\al,t),t)\da z_1(\al,t) d\al dt\\
&=K_5+K_6,
\end{align*}
and $K_5$ is treated as  $K_3$ (a term in $K_5$ is part of
$F(x,t)$ and another of $f(\al,t)$).  $K_6$ can be written in the
following manner
\begin{align}
\begin{split}\label{k6}
K_6&=-(\rho^2\!-\!\rho^1)\int_0^\mathrm{T}\!\!\int_{-\pi}^\pi [BR_1 BR_2+
\frac{\varpi^2}{4}\frac{\da z_1\da z_2}{|\da z|^4}]\dxu\lambda(z)\da z_1 d\al dt\\
&\quad-(\rho^2\!+\!\rho^1)\int_0^\mathrm{T}\!\!\int_{-\pi}^\pi [\frac{\varpi}{2}BR_1\frac{\da z_2}{|\da z|^2}
+\frac{\varpi}{2}BR_2\frac{\da z_1}{|\da z|^2}]\dxu\lambda(z)\da z_1 d\al dt.
\end{split}
\end{align}
Regarding  $J_5$ we have the splitting
\begin{align*}
J_5&=\int_0^\mathrm{T}\!\!\int_{\R^2}\rho (v_2\dxu v_1 +v_1\dxu v_2)\dxd\lambda dxdt\\
&\quad+\int_0^\mathrm{T}\!\!\int_{-\pi}^\pi\big(\rho^2(v^2_1v^2_2)(z(\al,t),t)
-\rho^1(v_1^1v_2^1)(z(\al,t),t)\big)\dxd\lambda(z(\al,t),t)\da z_2(\al,t) d\al dt\\
&=K_7+K_8.
\end{align*}
 $K_7$ again can be treated like  $K_3$, and we obtain for $K_8$ the
 following expression
\begin{align}
\begin{split}\label{k8}
K_8&=(\rho^2\!-\!\rho^1)\int_0^\mathrm{T}\!\!\int_{-\pi}^\pi [BR_1 BR_2+
\frac{\varpi^2}{4}\frac{\da z_1\da z_2}{|\da z|^4}]\dxd\lambda(z)\da z_2 d\al dt\\
&\quad +(\rho^2\!+\!\rho^1)\int_0^\mathrm{T}\!\!\int_{-\pi}^\pi [\frac{\varpi}{2}BR_1\frac{\da z_2}{|\da z|^2}
+\frac{\varpi}{2}BR_2\frac{\da z_1}{|\da z|^2}]\dxd\lambda(z)\da z_2 d\al dt.
\end{split}
\end{align}
Next for $J_6$
\begin{align*}
J_6&=\int_0^\mathrm{T}\!\!\int_{\R^2}2\rho v_2\dxd v_2\dxd\lambda dxdt\\
&\quad-\int_0^\mathrm{T}\!\!\int_{-\pi}^\pi\big(\rho^2(v^2_2(z(\al,t),t))^2
-\rho^1(v_2^1(z(\al,t),t))^2\big)
\dxd\lambda(z(\al,t),t)\da z_2(\al,t) d\al dt\\
&=K_9+K_{10},
\end{align*}
and for $K_9$ we proceed as before. Finally we have
\begin{align}
\begin{split}\label{k10}
K_{10}&=-(\rho^2\!-\!\rho^1)\int_0^\mathrm{T}\!\!\int_{-\pi}^\pi [(BR_2)^2+
\frac{\varpi^2}{4}\frac{(\da z_2)^2}{|\da z|^4}]\dxd\lambda(z)\da z_1 d\al dt\\
&\quad -(\rho^2\!+\!\rho^1)\int_0^\mathrm{T}\!\!\int_{-\pi}^\pi \varpi BR_2\frac{\da z_2}{|\da z|^2}\dxd\lambda(z)\da z_1 d\al dt.
\end{split}
\end{align}
Using equations (\ref{k4}--\ref{k10}) we get the following sum
$K_4+K_6+K_8+K_{10}=(\rho^2-\rho^1)L_1+(\rho^2+\rho^1)L_2$
where
$$
L_1=-\int_0^\mathrm{T}\!\!\int_{-\pi}^\pi BR(z,\varpi)(\al,t)\cdot
\grad\lambda(z(\al,t),t) BR(z,\varpi)(\al,t)\cdot\dpa z(\al,t)d\al dt,
$$and
$$
L_2=-\int_0^\mathrm{T}\!\!\int_{-\pi}^\pi \frac{\varpi(\al,t)}{2|\da z(\al,t)|^2}
BR(z,\varpi)(\al,t)\cdot\dpa z(\al,t) \da z(\al,t)\cdot
\grad\lambda(z(\al,t),t) d\al dt.
$$
An integration by parts in the variable $\al$ in $L_2$ gives the last term of $f(\al,t)$.
 Then  identity \eqref{vnormal} gives $K_1+(\rho^2\!-\!\rho^1)L_1=0$ and the stated
  formula for the Laplacian of $p$ is proved. \quad\eop\\
\newline

Identity \eqref{vnormal} allows us to choose the velocity of the
curve as follows:
\begin{align}\label{eev}
\begin{split}
z_t(\al,t)=BR(z,\varpi)(\al,t)+c(\al,t)\da z(\al,t),
\end{split}
\end{align}
where the scalar $c(\al,t)$ is given by
\begin{align}
\begin{split}\label{fc2}
c(\al)&=\frac{\al+\pi}{2\pi}\int_{-\pi}^\pi\frac{\partial_{\beta}
z(\beta)}{|\partial_{\beta}
z(\beta)|^2}\cdot \partial_{\beta} BR(z,\varpi)(\beta) d\beta-\int_{-\pi}^\al \frac{\partial_{\beta}
z(\beta)}{|\partial_{\beta} z(\beta)|^2}\cdot
\partial_{\beta} BR(z,\varpi)(\beta) d\beta,
\end{split}
\end{align}
and has been  taken in such a way that the length of the tangent vector only depends on
the variable $t$:
\begin{equation}\label{cancelacionextra}
|\da z(\al,t)|^2=A(t).
\end{equation}
Since $c(\al,t)$ has to be periodic, we obtain
\begin{align}
\begin{split}\label{AppA}
A'(t)=2\da z_t(\al,t)\cdot\da z(\al,t)=\frac{1}{\pi}\int_{-\pi}^\pi\da z(\al,t)\cdot \da BR(z,\varpi)(\al,t)
d\al.
\end{split}
\end{align}

Next we close the system giving the evolution equation for the
amplitude of the vorticity $\varpi(\al,t)$ by means of Bernoulli's
law. This fact allows us to satisfy \eqref{weakeuler} showing that
we have a weak solution. Using \eqref{vfc} for $x\neq z(\al,t)$ we
get $v(x,t)=\grad \phi(x,t)$ where
$$
\phi(x,t)=\frac{1}{2\pi}PV\int
\arctan\Big(\frac{x_2-z_2(\beta,t)}{x_1-z_1(\beta,t)}\Big)\varpi(\beta,t)d\beta.
$$
Let us  define
$$
\Pi(\al,t)=\phi^2(z(\al,t),t)-\phi^1(z(\al,t),t),
$$
where again $\phi^j(z(\alpha,t),t)$ denotes the limit obtained
approaching the boundary in the normal direction inside $\Omega^j$.
It is clear that
\begin{equation*}
\begin{array}{rll}
\da \Pi(\al,t)&=(\grad\phi^2(z(\al,t),t)-\grad\phi^1(z(\al,t),t))\cdot \da z(\al,t)&\\
&=(v^2(z(\al,t),t)-v^1(z(\al,t),t))\cdot \da z(\al,t)&=\varpi(\al,t),
\end{array}
\end{equation*}
and therefore
$$
\int_{-\pi}^\pi \varpi(\al,t) d\al=0.
$$
Now we observe that
\begin{align}\label{limitepotencial}
\begin{split}
\phi^2(z(\alpha,t),t)&=IT(z,\varpi)(\al,t)
+\frac{1}{2}\Pi(\al,t),\\
\phi^1(z(\alpha,t),t)&=IT(z,\varpi)(\al,t)-\frac{1}{2}\Pi(\al,t),
\end{split}
\end{align}
where
$$
IT(z,\varpi)(\al,t)=\frac{1}{2\pi}PV\int
\arctan\Big(\frac{z_2(\al,t)-z_2(\beta,t)}{z_1(\al,t)-z_1(\beta,t)}\Big)\varpi(\beta,t)d\beta.
$$

Using the Bernoulli's law in \eqref{euler}, inside each domain, we
have
$$
\rho(\phi_t(x,t)+\frac{1}{2}|v(x,t)|^2+\mathrm{g}x_2)+p(x,t)=0.
$$
Next we take  limits to get
$$
\rho^j(\phi^j_t(z(\al,t),t)+\frac{1}{2}|v^j(z(\al,t),t)|^2+\mathrm{g}z_2(\al,t))+p^j(z(\al,t),t)=0,
$$
and since $ p^1(z(\al,t),t)=p^2(z(\al,t),t),$ we obtain
\begin{equation}\label{limiteB}
[\rho\phi_t](\al,t)+\frac{\rho^2}{2}|v^2(z(\al,t),t)|^2-\frac{\rho^1}{2}|v^1(z(\al,t),t)|^2
+(\rho^2-\rho^1)\mathrm{g}z_2(\al,t)=0,
\end{equation}
where we have introduced the following notation: $$
[\rho\phi_t](\al,t)=\rho^2\phi^2_t(z(\al,t),t)-\rho^1\phi^1_t(z(\al,t),t).$$
Then it is clear that
$\phi^j_t(z(\al,t),t)=\partial_t(\phi^j(z(\al,t),t))-z_t(\al,t)\cdot\grad\phi^j(z(\al,t),t),$
and using \eqref{limitepotencial} we find that
\begin{align*}
[\rho\phi_t]=\frac{\rho^2+\rho^1}{2}\Pi_t+(\rho^2-\rho^1)\partial_t(IT(z,\varpi))-z_t\cdot(\rho^2 v^2(z,t)-\rho^1
v^1(z,t)).
\end{align*}
Introducing equations \eqref{limitevelocidad} and \eqref{eev} into
\eqref{limiteB} we get
\begin{align}
\begin{split}\label{eep}
\Pi_t(\al,t)&=-2A_\rho\partial_t(IT(z,\varpi))(\al,t) +c(\al,t)\varpi(\al,t)+A_\rho|BR(z,\varpi)(\al,t)|^2\\
&\quad +2A_\rho c(\al,t)BR(z,\varpi)(\al,t)\cdot \da
z(\al,t)-A_\rho\frac{|\varpi(\al,t)|^2}{4|\da
z(\al,t)|^2}-2A_\rho\mathrm{g}z_2(\al,t).
\end{split}
\end{align}
Since the equality
\begin{align*}
\da\partial_t(IT(z,\varpi))=\dpt(BR(z,\varpi)\cdot\da z)&=\dpt BR(z,\varpi)\cdot\da z+BR(z,\varpi)\cdot \da BR(z,\varpi)\\
&\quad +cBR(z,\varpi)\cdot \da^2 z+\da cBR(z,\varpi)\cdot \da z
\end{align*}
can be proved easily, we can take then  a derivative in \eqref{eep} and use the above identity
to find the desired formula for $\varpi$:
\begin{align}
\begin{split}\label{fw2}
\varpi_t(\al,t)&=-2A_\rho \partial_t BR(z,\varpi)(\al,t)\cdot
\da z(\al,t)-A_\rho\da( \frac{|\varpi|^2}{4|\da z|^2})(\al,t)+\da (c\, \varpi)(\al,t)\\
&\quad +2A_{\rho}c(\al,t)\da BR(z,\varpi)(\al,t)\cdot\da
z(\al,t)-2A_\rho\mathrm{g}\da z_2(\al,t).
\end{split}
\end{align}\\

Our next step will be to get the formula for the difference of the
gradients of the pressure in the normal direction:
\begin{equation}\label{R-T}
\sigma(\al,t)=-(\nabla p^2(z(\alpha,t),t) - \nabla
p^1(z(\alpha,t),t))\cdot\partial^{\perp}_{\alpha}z(\alpha,t),
\end{equation}
which we shall find in the singular terms of the evolution equation.

 We will consider the case $\rho_1=0$, which gives $-\nabla p(x,t)=0$ inside
$\Omega^1(t)$ and therefore $\nabla p^1(z(\alpha,t),t)=0$. Let us
define the Lagrangian coordinates for the free boundary with the
velocity $v^2$
\begin{align*}
\begin{split}
Z_t(\gamma,t)&=v^2(Z(\gamma,t),t))\\
Z(\gamma,0)&=z_0(\gamma).
\end{split}
\end{align*}
We have  two different parameterizations for the same curve
$Z(\gamma,t)=z(\al(\gamma,t),t)$ and also two equations for its
velocity, namely
\begin{align*}
Z_t(\gamma,t)&=z_t(\al,t)+\al_t(\gamma,t)\da z(\al,t)\\
&=BR(z,\varpi)(\al,t)+c(\al,t)\da z(\al,t)+\al_t(\gamma,t)\da
z(\al,t)
\end{align*}
and another one given by the limit
\begin{align}\label{Zt}
Z_t(\gamma,t)&=BR(z,\varpi)(\al,t)+\frac12\frac{\varpi(\al,t)}{|\da
z(\al,t)|^2}\da z(\al,t).
\end{align}
The dot product with the tangential vector gives
\begin{equation*}
\al_t(\gamma,t)=\frac12\frac{\varpi(\al,t)}{|\da
z(\al,t)|^2}-c(\al)=\frac{\varphi(\al,t)}{|\da z(\al,t)|}.
\end{equation*}
And taking a  time derivative in \eqref{Zt} yields
\begin{align*}
Z_{tt}(\gamma,t)\cdot \dpa z(\al,t)&= (\partial_t BR(z,\varpi)(\al,t)+\al_t(\gamma,t)\da BR(z,\varpi)(\al,t))\cdot \dpa z(\al,t)\\
&\quad+\frac12\frac{\varpi(\al,t)}{|\da z(\al,t)|^2}(\da
z_t(\al,t)+\al_t(\gamma,t)\da^2 z(\al,t))\cdot \dpa z(\al,t)
\end{align*}
Therefore
\begin{align}
\begin{split}\label{fsigma}
\frac{\sigma(\al,t)}{\rho^2}&=(\partial_t BR(z,\varpi)(\al,t)+\frac{\varphi(\al,t)}{|\da z(\al,t)|}\da BR(z,\varpi)(\al,t))\cdot \dpa z(\al,t)\\
&\quad+\frac12\frac{\varpi(\al,t)}{|\da z(\al,t)|^2}(\da
z_t(\al,t)+\frac{\varphi(\al,t)}{|\da z(\al,t)|}\da^2 z(\al,t))\cdot
\dpa z(\al,t)+\mathrm{g}\da z_1(\al,t).
\end{split}
\end{align}

\begin{rem}
Let us consider $\rho_2$ and $\rho_1$ to be now arbitrary
densities, then using the lagrangian coordinates for the free
boundary of the fluid in $\Omega^1(t)$
\begin{align*}
\begin{split}
Z'_t(\gamma,t)&=v^1(Z'(\gamma,t),t))\\
Z'(\gamma,0)&=z_0(\gamma),
\end{split}
\end{align*}
it is easy to check that
\begin{align*}
\begin{split}
\frac{\sigma(\al,t)}{\rho^2+\rho^1}&=A_\rho(\partial_t BR(z,\varpi)(\al,t)+
\frac{|\varpi(\al,t)|^2}{4|\da z(\al,t)|^4}\da^2 z(\al,t))\cdot \dpa z(\al,t)\\
&\quad+(\frac{\varpi(\al,t)}{|\da z(\al,t)|^2}-A_\rho c(\al,t))\da
BR(z,\varpi)(\al,t)\cdot\dpa z(\al,t)+\mathrm{g}A_\rho \da
z_1(\al,t).
\end{split}
\end{align*}
\end{rem}

%%%%%%%%%%%%%%%%%%%%%%%%%%%%%%%%%%%%%%%%%%%%%%%%%%%%%%%%%%%%%%%%%%%%%%%%%%%%%%%%%%%
%%%%%%%%%%%%%%%%%%%%%%%%%%%%%%%%%%%%%%%%%%%%%%%%%%%%%%%%%%%%%%%%%%%%%%%%%%%%%%%%%%%

\section{The evolution equation in terms of $\varphi(\al,t)$}

%%%%%%%%%%%%%%%%%%%%%%%%%%%%%%%%%%%%%%%%%%%%%%%%%%%%%%%%%%%%%%%%%%%%%%%%%%%%%%%%%%%
%%%%%%%%%%%%%%%%%%%%%%%%%%%%%%%%%%%%%%%%%%%%%%%%%%%%%%%%%%%%%%%%%%%%%%%%%%%%%%%%%%%

We will consider $\rho_1=0$ and therefore $A_\rho=1$. Using
\eqref{fw2} we can write

\begin{align}
\begin{split}\label{fw3}
\varpi_t(\al,t)&=-2\partial_t BR(z,\varpi)(\al,t)\cdot
\da z(\al,t)-\da( \frac{|\varpi|^2}{4|\da z|^2})(\al,t) +\da (c\, \varpi)(\al,t)\\
&\quad +2c(\al,t)\da BR(z,\varpi)(\al,t)\cdot\da
z(\al,t)-2\mathrm{g}\da z_2(\al,t),
\end{split}
\end{align}

%\begin{align}
%\begin{split}\label{fw4}
%\varpi_t(\al,t)+T(\varpi_t)(\al,t)&=R(\al,t)-\da( \frac{|\varpi|^2}{4|\da z|^2})(\al,t)+\da (c\, \varpi)(\al,t)\\
%&\quad +2c(\al,t)\da BR(z,\varpi)(\al,t)\cdot\da
%z(\al,t)+2\mathrm{g}\da z_2(\al,t).
%\end{split}
%\end{align}
%Therefore
%\begin{align*}
%\varpi_t(\al,t)&=(I+T)^{-1}(f)(\al,t),
%\end{align*}
%where $f=f(z,\varpi)$. From \cite{BMO} we know that the operator
%is invertible, since in the proof we will use energy estimates, we
%will need to estimate $\varpi_t$ in terms of $z,\varpi,\varphi$
%\eqref{nswt}. So we will need Theorem \ref{oil2l2}.

In the case $A_\rho=0$ the expression \eqref{fw2} yields
$$
\varpi_t(\al,t)=\da (c\, \varpi)(\al,t),
$$
that is, we are obtain the vortex sheet problem for which the
Kelvin-Helmholtz instability arises \cite{Ebin2} \cite{ACG}. For $A_\rho=1$
this term again appears in the evolution equation, and in order to
absorb it we shall make use of the parameter $\varphi(\al,t)$
\cite{BHL} \cite{AM}. The fact that $|\da z(\al,t)|^2=A(t)$ yields
$$
2A(t)\da c=\frac{1}{\pi}\int_{-\pi}^\pi \!\da z(\al,t)\cdot \da
BR(z,\varpi)(\al,t) d\al-2\da z\cdot \partial_{\al}
BR(z,\varpi)
$$
and therefore
$$
2c\,\da z\cdot \partial_{\al} BR(z,\varpi)=-\da
(A\, c^2)+\frac{c}{\pi}\int_{-\pi}^\pi \!\da z(\al,t)\cdot \da
BR(z,\varpi)(\al,t) d\al.
$$
Substituting  the  formula above in \eqref{fw3} we find
\begin{align}
\begin{split}\label{nose}
\varpi_t&=-2\partial_t BR(z,\varpi)\cdot
\da z-\da(\varphi^2)+\frac{c}{\pi}\int_{-\pi}^\pi \!\da z(\al,t)\cdot \da
BR(z,\varpi)(\al,t) d\al-2\mathrm{g}\da z_2,
\end{split}
\end{align}
for $\varphi$ given by \eqref{fphi}. From that identity we have
\begin{equation}\label{phit}
\varphi_t(\al,t)=\frac{\varpi_t(\al,t)}{2|\da
z(\al,t)|}-\frac{\varpi(\al,t)}{2|\da z(\al,t)|^3}\da z(\al,t)\cdot\da z_t(\al,t)-\partial_t(c|\da
z|)(\al,t)
\end{equation}
which together with \eqref{nose} and \eqref{AppA} yields
\begin{align*}
\begin{split}
\varphi_t=-\partial_t BR(z,\varpi)\cdot
\frac{\da z}{|\da z|}&-\frac{\da(\varphi^2)}{2|\da z|}+c\,\frac{1}{2\pi}\int_{-\pi}^\pi \!\frac{\da z(\al,t)}{|\da
z(\al,t)|}\cdot \da BR(z,\varpi)(\al,t) d\al-\mathrm{g}\frac{\da z_2}{|\da z|}\\
&-\frac{\varpi}{2|\da z|^2}\frac{1}{2\pi}\int_{-\pi}^\pi \!\frac{\da z(\al,t)}{|\da
z(\al,t)|}\cdot \da BR(z,\varpi)(\al,t)
d\al-\partial_t(c|\da z|),
\end{split}
\end{align*}
that is
\begin{align}
\begin{split}\label{ephi}
\varphi_t&=-\frac{\da(\varphi^2)}{2|\da
z|}-B(t) \,\varphi-\partial_t
BR(z,\varpi)\cdot
\frac{\da z}{|\da z|}-\mathrm{g} \frac{\da z_2}{|\da
z|}-\partial_t(c|\da z|).
\end{split}
\end{align}
where
$$
B(t)=\frac{1}{2\pi}\int_{-\pi}^\pi \!\frac{\da z(\al,t)}{|\da
z(\al,t)|^2}\cdot \da BR(z,\varpi)(\al,t) d\al.
$$
It is easy to check in the equation above that the singular term
$\da (c\, \varpi)$ takes part of the transport term
$\da(\varphi^2)$.

Now let us remember that the evolution equation for the quantity
$\Pi(\al,t)$  was discovered using the continuity of the pressure on
$z(\al,t)$ (Proposition \ref{presiones}). Analogously the evolution
equation for $\da\Pi(\al,t)=\varpi(\al,t)$ can be obtained
throughout the following identity:
$$
-(\grad p^2(z(\al,t),t)-\grad p^1(z(\al,t),t))\cdot\da z(\al,t)=0.
$$
(Observe nevertheless that the Rayleigh-Taylor condition refers the
jump of the pressure in the  normal direction \eqref{R-T}.)

 With the help of
property \eqref{cancelacionextra} we find that
$$
\da^2z(\al,t)\cdot\da z(\al,t)=0,
$$
and therefore
$$
\da^2 z(\al,t)=\frac{\da^2z(\al,t)\cdot\dpa z(\al,t)}{|\da z(\al,t)|^2}\dpa z(\al,t).
$$
In the above formula we get the normal direction in the second
derivative of $z$. Using this fact in \eqref{ephi} we obtain
\begin{align*}
\da\varphi_t&=-\frac{\da^2(\varphi^2)}{2|\da
z|}-B(t) \,\da\varphi -(\partial_t
BR(z,\varpi)\cdot\dpa z+\mathrm{g}\da z_1)  \frac{\da^2z\cdot\dpa z}{|\da z|^3}\\
&\quad -\partial_t\Big(\frac{1}{2\pi}\int_{-\pi}^\pi\!\frac{\da z(\al,t)}{|\da
z(\al,t)|}\cdot \da BR(z,\varpi)(\al,t) d\al\Big)+\partial_t(\frac{\da z}{|\da z|})\cdot\da BR(z,\varpi).
\end{align*}
In $\partial_t(\frac{\da z}{|\da z|})$  the perpendicular direction
also appears, so that completing the formula for $\sigma$
\eqref{fsigma} we get
\begin{align*}
\da\varphi_t&=-\frac{\da^2(\varphi^2)}{2|\da
z|}-B(t) \,\da\varphi -\frac{\sigma}{\rho^2}\frac{\da^2z\cdot\dpa z}{|\da z|^3}
-\partial_t(|\da z|B)(t)+\da BR(z,\varpi)\cdot\dpa z\frac{\da z_t\cdot\dpa z}{|\da z|^3}\\
&\quad +\Big(\frac{\varphi}{|\da z|}\da BR(z,\varpi)\cdot\dpa z+\frac12\frac{\varpi}{|\da z|^2}(\da z_t\cdot\dpa z+\frac{\varphi}{|\da z|}\da^2 z\cdot\dpa z)\Big)\frac{\da^2z\cdot\dpa z}{|\da z|^3}
\end{align*}
and therefore
\begin{align*}
\da\varphi_t&=-\frac{\da^2(\varphi^2)}{2|\da
z|}-B(t) \,\da\varphi -\frac{\sigma}{\rho^2}\frac{\da^2z\cdot\dpa z}{|\da z|^3}-\partial_t(|\da z|B)(t)\\
&\quad +\frac{1}{|\da z|^3}\Big(\da BR(z,\varpi)\cdot\dpa z+\frac{\varpi}{2|\da z|^2}\da^2z\cdot\dpa z\Big)\Big(\da z_t\cdot\dpa z+\frac{\varphi}{|\da z|}\da^2 z\cdot\dpa z\Big).
\end{align*}
Finally after a straightforward calculation we obtain the following:
\begin{align}
\begin{split}\label{eephi}
\da\varphi_t&\!=-\frac{\da^2(\varphi^2)}{2|\da
z|}-B\,\da\varphi-\frac{\sigma}{\rho^2}\frac{\da^2z\cdot\dpa z}{|\da z|^3}-\partial_t(|\da z|B)\\
&\quad+\frac{1}{|\da z|^3}(\da BR(z,\varpi)\cdot \dpa
z+\frac{\varpi}{2|\da z|^2}\da^2 z\cdot \dpa z)^2.
\end{split}
\end{align}

%%%%%%%%%%%%%%%%%%%%%%%%%%%%%%%%%%%%%%%%%%%%%%%%%%%%%%%%%%%%%%%%%%%%%%%%%%%%%%%%%%%%%%%%%%%%%%%%%%%%%%%%%%%%%%%%
%%%%%%%%%%%%%%%%%%%%%%%%%%%%%%%%%%%%%%%%%%%%%%%%%%%%%%%%%%%%%%%%%%%%%%%%%%%%%%%%%%%%%%%%%%%%%%%%%%%%%%%%%%%%%%%%

\section{The basic operator}

%%%%%%%%%%%%%%%%%%%%%%%%%%%%%%%%%%%%%%%%%%%%%%%%%%%%%%%%%%%%%%%%%%%%%%%%%%%%%%%%%%%%%%%%%%%%%%%%%%%%%%%%%%%%%%%%
%%%%%%%%%%%%%%%%%%%%%%%%%%%%%%%%%%%%%%%%%%%%%%%%%%%%%%%%%%%%%%%%%%%%%%%%%%%%%%%%%%%%%%%%%%%%%%%%%%%%%%%%%%%%%%%%

Let the operator $T$ be defined by the formula

\begin{equation}
T(u)(\al)=2BR(z,u)(\al)\cdot \da z(\al). \label{operatorT}
\end{equation}

\begin{lemma}
Suppose that $\|\F(z)\|_{L^{\infty}}< \infty$ \eqref{df} and $z\in
C^{2,\delta}$ with $0<\delta<1/2$. Then $T:L^2\rightarrow H^1$ and
\begin{equation}\label{ol2h1}
\|T\|_{L^2\rightarrow H^1}\leq \|\F(z)\|^4_{L^{\infty}}\|z\|^4_{C^{2,\delta}}.
\end{equation}
\end{lemma}

Proof: Here we shall show the argument in the case of a closed
curve. The other case was treated in \cite{DY2}.

Since the formula \eqref{fibr} yields
$$T(u)(\al)=\frac1\pi\da \int_{-\pi}^\pi u(\beta)\arctan \Big(\frac{z_2(\al)-z_2(\beta)}{z_1(\al)-z_1(\beta)}\Big)d\beta,$$
we have
$$\int_{-\pi}^\pi T(u)(\al)d\al=0, $$
which implies $\|T(u)\|_{L^2}\leq \|\da T(u)\|_{L^2}$.

Let us write first: $$\da T(u)=2BR(z,u)(\al)\cdot \da^2 z(\al)+2\da
z(\al) \cdot \da BR(z,u)(\al)=I_1+I_2.$$ For $I_1$ we have the
expression
\begin{align*}
\begin{split}
I_1&=2(BR(z,u)(\al)-\frac{\dpa z(\al)}{|\da z(\al)|^2} H(u)(\al))\cdot \da^2 z(\al)+2H(u)(\al)\frac{\dpa z(\al)\cdot \da^2 z(\al)}{|\da z(\al)|^2}\\
&=J_1+J_2,
\end{split}
\end{align*}
where $H(u)$ is the (periodic) Hilbert transform of the function u.

Then
\begin{align*}
\begin{split}
J_1&=\frac{1}{\pi}\da^2
z(\al)\cdot\int_{-\pi}^\pi u(\al-\beta)[\frac{(z(\al)-z(\al-\beta))^\bot}{|z(\al)-z(\al-\beta)|^2}-\frac{\dpa z(\al)}{2|\da z(\al)|^2\tan(\beta/2)}] d\beta.
\end{split}
\end{align*}
Let us define
\begin{equation}\label{fb}
C_1(\al,\beta)=\frac{(z(\al)-z(\al-\beta))^\bot}{|z(\al)-z(\al-\beta)|^2}-\frac{\dpa z(\al)}{2|\da z(\al)|^2\tan(\beta/2)},
\end{equation} then we shall show that $\|C_1\|_{L^\infty}\leq C\|\F(z)\|^2_{L^{\infty}}\|z\|^2_{C^2}$ and therefore $J_1\leq C\|\F(z)\|2_{L^{\infty}}\|z\|^3_{C^2}\|u\|_{L^2}$. Since the
estimate  $J_2\leq C\|\F(z)\|_{L^{\infty}}\|z\|_{C^2}|H(u)(\al)|$ is immediate,
we finally get
\begin{align}\label{oti1}
\begin{split}
|I_1|\leq
C\|\F(z)\|^2_{L^{\infty}}\|z\|^3_{C^2}(\|u\|_{L^2}+|H(u)(\al)|).
\end{split}
\end{align}
Next we split $C_1=D_1+D_2+D_3$ where
$$
D_1=\frac{(z(\al)-z(\al-\beta)-\da z(\al)\beta)^\bot}{|z(\al)-z(\al-\beta)|^2},\quad D_2=\dpa z(\al)[\frac{\beta}{|z(\al)-z(\al-\beta)|^2}-\frac{1}{|\da z(\al)|^2\beta}],
$$
and
$$
D_3=\frac{\dpa z(\al)}{|\da z(\al)|^2}[\frac1\beta-\frac{1}{2\tan(\beta/2)}].
$$
The inequality \begin{equation}\label{ncll}|z(\al)-z(\al-\beta)-\da
z(\al)\beta|\leq \|z\|_{C^2}|\beta|^2\end{equation} yields easily
$|D_1|\leq \|z\|_{C^2}\|\F(z)\|^2_{L^{\infty}}$.

Then we can rewrite $D_2$ as follows:
$$
D_2=\dpa z(\al)[\frac{(\da z(\al)\beta-(z(\al)-z(\al-\beta)))\cdot (\da z(\al)\beta+(z(\al)-z(\al-\beta)))}{|z(\al)-z(\al-\beta)|^2|\da z(\al)|^2\beta}],
$$
and, in particular, we have
$$
|D_2|\leq \frac{|\da z(\al)\beta-(z(\al)-z(\al-\beta))|(|\da z(\al)\beta|+|z(\al)-z(\al-\beta)|)}{|z(\al)-z(\al-\beta)|^2|\da z(\al)||\beta|}.
$$
Using \eqref{ncll} we find that $|D_2|\leq 2
\|z\|_{C^2}\|\F(z)\|^2_{L^{\infty}}$.

Next let us observe that $[-\pi,\pi]$ gives $|D_3|\leq
C\|\F(z)\|_{L^{\infty}}$.

 The identity $\da z(\al)\cdot\dpa z(\al)=0$ allows us to write
$I_2$ as follows: $$
I_2=-\frac{2}{\pi}\int_{-\pi}^\pi u(\beta)\frac{(z(\al)-z(\beta))^\bot\cdot\da
z(\al)(z(\al)-z(\beta))\cdot \da z(\al)}{|z(\al)-z(\beta)|^4}d\beta.
$$
and therefore $$
I_2=-\frac{2}{\pi}\int_{-\pi}^\pi u(\al-\beta)\frac{(z(\al)-z(\al-\beta)-\da z(\al)\beta)^\bot\cdot\da
z(\al)(z(\al)-z(\al-\beta))\cdot \da z(\al)}{|z(\al)-z(\al-\beta)|^4}d\beta.
$$
Next we take $I_2=J_3+J_4+J_5+J_6+J_7$ where
$$
J_3=-\frac{2}{\pi}\int_{-\pi}^\pi u(\al-\beta)\frac{(E(\al,\beta))^\bot\cdot\da
z(\al)(z(\al)-z(\al-\beta))\cdot \da z(\al)}{|z(\al)-z(\al-\beta)|^4}d\beta,
$$
$$
J_4=-(\da^2z(\al))^\bot\cdot\da
z(\al)\frac{1}{\pi}\int_{-\pi}^\pi u(\al-\beta)\frac{\beta^2((z(\al)-z(\al-\beta)-\da z(\al)\beta)\cdot \da z(\al)}{|z(\al)-z(\al-\beta)|^4}d\beta,
$$
$$
J_5=-(\da^2z(\al))^\bot\cdot\da
z(\al)|\da z(\al)|^2\frac{1}{\pi}\int_{-\pi}^\pi u(\al-\beta)[\frac{\beta^3}{|z(\al)-z(\al-\beta)|^4}-\frac{1}{|\da z(\al)|^4\beta}]d\beta,
$$
$$
J_6+J_7=-\frac{(\da^2z(\al))^\bot\cdot\da
z(\al)}{|\da z(\al)|^2}\Big(\frac{1}{\pi}\int_{-\pi}^\pi u(\al-\beta)[\frac{1}{\beta}-\frac{1}{2\tan(\beta/2)}]d\beta+H(u)(\al)\Big),
$$ and $E(\al,\beta)=z(\al)-z(\al-\beta)-\da z(\al)\beta-\frac12\da^2z(\al)\beta^2$.
Using the  bound
\begin{equation}\label{holder}
|E(\al,\beta)|\leq \frac12\|z\|_{C^{2,\delta}}|\beta|^{2+\delta},
\end{equation}
one get easily that $$|J_3|\leq
C\|z\|^3_{C^{2,\delta}}\|\F(z)\|^3_{L^{\infty}}\int_{-\pi}^\pi|u(\beta)||\beta|^{\delta-1}.$$
Then reasoning as before the  inequality \eqref{ncll} gives as
$|J_4|\leq C\|z\|^4_{C^2}\|\F(z)\|^4_{L^{\infty}}\|u\|_{L^2}$.
Regarding $D_2$, we have $|J_5|\leq C
\|z\|^4_{C^2}\|\F(z)\|^3_{L^\infty}\|u\|_{L^2}$, and it is easy to
get $|J_6|\leq C\|z\|_{C^2}\|\F(z)\|_{L^\infty}\|u\|_{L^2}$. Finally
we have
$$
|I_2|\leq
C\|\F(z)\|^4_{L^{\infty}}\|z\|^4_{C^{2,\delta}}(\|u\|_{L^2}+|H(u)(\al)|+\int_{-\pi}^\pi
|\beta|^{\delta-1}|u(\al-\beta)|d\beta).
$$
This last inequality together with \eqref{oti1} gives us
$$
|\da T(u)(\al)|\leq
C\|\F(z)\|^4_{L^{\infty}}\|z\|^4_{C^{2,\delta}}(\|u\|_{L^2}+|H(u)(\al)|+\int_{-\pi}^\pi
|\beta|^{\delta-1}|u(\al-\beta)|d\beta).
$$
To finish we use the  $L^2$ boundedness of $H$ and  Minkowski's
inequality to obtain the estimate
$$
\|\da T(u)\|_{L^2}\leq
C\|\F(z)\|^4_{L^{\infty}}\|z\|^4_{C^{2,\delta}}\|u\|_{L^2},
$$
q.e.d.

%%%%%%%%%%%%%%%%%%%%%%%%%%%%%%%%%%%%%%%%%%%%%%%%%%%%%%%%%%%%%%%%%%%%%%%%%%%%%%%%%%%%%%%%%%%%%%%%%%%%%%%%%%%%%%%%%%%%%%%%%%
%%%%%%%%%%%%%%%%%%%%%%%%%%%%%%%%%%%%%%%%%%%%%%%%%%%%%%%%%%%%%%%%%%%%%%%%%%%%%%%%%%%%%%%%%%%%%%%%%%%%%%%%%%%%%%%%%%%%%%%%%%

\section{Estimates on the inverse operator $(I+T)^{-1}$. }

%%%%%%%%%%%%%%%%%%%%%%%%%%%%%%%%%%%%%%%%%%%%%%%%%%%%%%%%%%%%%%%%%%%%%%%%%%%%%%%%%%%%%%%%%%%%%%%%%%%%%%%%%%%%%%%%%%%%%%%%%%%
%%%%%%%%%%%%%%%%%%%%%%%%%%%%%%%%%%%%%%%%%%%%%%%%%%%%%%%%%%%%%%%%%%%%%%%%%%%%%%%%%%%%%%%%%%%%%%%%%%%%%%%%%%%%%%%%%%%%%%%%%%

As it was shown in reference \cite{DY2}, under our hypothesis about
the curve $z$, $T(u)=2BR(z,u)(\al)\cdot \da z(\al)$ defines a
compact operator in  Sobolev spaces. Its adjoint $T^*$ is given as
the real part of the Cauchy integral and it does not has real
eigenvalues $\lambda$ such that $|\lambda|\geq1$ \cite{BMO}.
Therefore the existence of the bounded operator $(I+T)^{-1}$ follows
from the standard theory.

Let $F(z)$ given by $$F(z)= \frac{1}{2\pi i} \int \frac{u(\beta)\da
z(\beta)}{z - z(\beta)}d\beta,$$ and $f(z)= Re(F(z))$, which can be
considered either in the periodic setting, where we have two
periodic domains $\Omega^1$, $\Omega^2$ (see ref. \cite{DY2}), or in
the bounded domain case ($\Omega^2$ bounded). In both situations
$F(z)$ can be evaluated in the interior of both domains, and $T^*$
appears when we take limits approaching the boundary from the
interior of each $\Omega^j$: $z=z(\al)+\ep\dpa z(\al)$,
$\epsilon\rightarrow 0$, ($\epsilon >0, \Omega^1; \epsilon <0,
\Omega^2$):
\begin{equation*}\label{fec}
f(z(\al))=T^*(u) - \sign(\ep)u(\al).
\end{equation*}

The periodic case was treated in ref.\cite{DY2} (proposition 4.2).
Therefore we shall consider here the bounded domain case.

Let $\mathcal{H}^j$ denote the Hilbert transform associated to
$\Omega^j$, we have:
\begin{align*}
\begin{split}
(\mathcal{H}^j)^2 = -I,\\
F^1 = F/\Omega^1 = f^1 + i g^1,\\
F^2 = F/\Omega^2 = f^2 + i g^2,\\
f^1/\partial \Omega = T^*u - u,\\
f^2/\partial \Omega = T^*u + u,\\
g^1/\partial \Omega =g^2/\partial \Omega = \mathcal{G}(u),\\
u - T^*u = \mathcal{H}^1(\mathcal{G}(u)),\\
u + T^*u = \mathcal{H}^2(\mathcal{G}(u)).
\end{split}
\end{align*}

\begin{thm}\label{oil2l2}
The norm of the operator $(I+T)^{-1}$ from $L^2$ to $L^2$ is bounded from above by
$\exp(C|||z|||^p)$ with
$|||z|||=\|z\|_{H^3}+\|\F(z)\|_{L^\infty}$, for some universal
constants $C$ and $p$.
\end{thm}
Proof: As in Proposition 4.2 (ref. \cite{DY2}) the proof follows from
the estimate
\begin{equation*}
||\mathcal{H}^j||_{L^2(\partial\Omega^j)}\leq \exp(C|||z|||^p)
\end{equation*}
Let $\phi$ be a conformal mapping of $\Omega^2$ into the unit disc
$D$ such that $\phi(z_0) = 0$ where $z_0$ satisfies $dist(z_0,
\partial\Omega^1)>>\frac{1}{|||z|||}$, then
\begin{eqnarray*}
\mathcal{H}^2f= H(f\circ\phi^{-1})\circ\phi
\end{eqnarray*}
where $H$ is the Hilbert transform in the unit disc D. Since
$\partial\Omega^2$ is smooth enough ($C^{2,\alpha}$) we know from
general theory that $\phi$ and $\phi'$ have continuous extensions
to $\partial\Omega^2$ and our problem is reduced to obtain a
weighted estimate for the Hilbert transform $H$ in $\partial D$
with respect to the weight $w(\tau) = |(\phi^{-1})'(\tau)|$,
$|\tau|=1$. But that is a consequence of the inequality
\begin{eqnarray*}
e^{-C|||z|||^p}\leq \frac{w(\tau_1)}{w(\tau_2)}\leq  e^{C|||z|||^p}
\end{eqnarray*}
for arbitrary $\tau_j$, $|\tau_j|=1$.

Following Riemann let us write $\phi(z) = (z - z_0)e^{R(z) + i
S(z)}$ where the real harmonic function $R(z)$ is the solution of
the following Dirichlet's problem:
\begin{eqnarray*}
\Delta R = 0 \quad \text{in} \quad \Omega^2\\
R(z)= -log|z - z_0|, \quad z\in \partial\Omega^2.
\end{eqnarray*}
Since $\Omega^2$ is a regular domain whose boundary has tangent
balls of radius $\frac{1}{|||z|||}$ contained in $\Omega^2$, it
follows from the standard theory that $|\nabla R|_{L^{\infty}}<<
|||z|||log(|||z|||)$. This estimate also holds for the conjugate
harmonic functions $S(z)$ implying $|\phi'(\tau)|<<
|||z|||log(|||z|||)$, $\tau\in\partial\Omega^2$.

Given $\tau_0\in \partial\Omega^2$ the arc $\gamma =
\{\tau\in\partial\Omega^2 : dist(\tau,\tau_0) <
\frac{1}{|||z|||log(|||z|||)}\}$ is then mapped by $\phi$ into the
semicircle $\phi(\gamma)=\{z\in\partial D : dist(z,\phi(\tau_0))\leq
\sqrt{2}\}$.

Let us consider the Cayley transform $C_{\phi(\tau_0)}: D
\rightarrow \mathbb{R}^2_+$
\begin{eqnarray*}
C_{\phi(\tau_0)}(z) = -\frac{1 - \overline{\phi(\tau_0)}\cdot z}{1 +
\overline{\phi(\tau_0)}\cdot z}
\end{eqnarray*}
verifying that
\begin{align*}
\begin{split}
V= Im(C_{\phi(\tau_0)}\circ\phi) \geq 0 \quad \text{in} \quad \Omega^2,\\
V/\partial\Omega^2 = 0,\\
w(\gamma) = Re(C_{\phi(\tau_0)}\circ\phi)(\gamma) \subset [-1, +1],\\
w(\tau_0) = 0.
\end{split}
\end{align*}

Applying Hopf's maximum principle to the non-negative harmonic
function $V$ in a disc of radius $1/|||z|||$ tangent to
$\partial\Omega^2$ in $\tau$, we get an estimate for the normal
derivative of $V$ at $\tau$ i.e. for $||\nabla V(\tau)||$ (since
$\partial\Omega^2$ is a level set of $V$), namely:
\begin{eqnarray*}
|\frac{\partial V}{\partial\nu}(\tau)| >> \frac{1}{|||z|||}V(\tau^*)
\end{eqnarray*}
where $\tau^*$ is the center of the disc.

To get an upper bound we may use the Poisson's kernel representation
of $V$ in a $C^{2,\alpha}$-domain $\widetilde{\Omega}$ contained in
$\Omega^2$ whose boundary consists of $\gamma$ and its parallel arc
$\gamma^*$  at distance $1/|||z|||$, together with two
"vertical" connecting arcs chosen in such a way that the
$C^{2,\alpha}$-norm of $\partial\widetilde{\Omega}$ is controlled by
$|||z|||$. Since $V/\partial \Omega^2\equiv 0 $ we obtain the
estimate:
\begin{eqnarray*}
|\frac{\partial V}{\partial\nu}(\tau)| << |||z|||log(|||z|||)
sup_{\tau\in\widetilde{\Omega}}V(\tau)
\end{eqnarray*}
for $$\tau\in \frac12\gamma = \{\tau\in\partial\Omega^2,
dist(\tau,\tau_0)\leq \frac{1}{2C|||z|||log(|||z|||)}\}.$$

We are then in condition to invoke Dahlberg's Harnack inequality up
to the boundary \cite{Dahlberg} to conclude that
\begin{eqnarray*}
|\frac{\partial V}{\partial\nu}(\tau)| << |||z|||log(|||z|||)
V(\tau^*), \quad \tau\in\frac12\gamma.
\end{eqnarray*}

Next we use the standard Harnack's inequality in the parallel curve
$\gamma^*$ to conclude that
\begin{eqnarray*}
\frac{||\nabla V(\tau_1)||}{||\nabla V(\tau_2)||} << |||z|||^2
log(|||z|||)
\end{eqnarray*}
for any two points $\tau_1,\tau_2 \in\frac12\gamma$.

But since $\frac12\leq|C'(\phi(\tau))|\leq 2$, $\tau\in \gamma$, we
get the bound
\begin{eqnarray*}
|\frac{\phi'(\tau_1)}{\phi'(\tau_2)}| << |||z|||^2 log(|||z|||)\quad
\tau_1,\tau_2\in\frac12\gamma.
\end{eqnarray*}

Let us observe now that the length of $\partial\Omega^2$ is
controlled by $|||z|||$ giving us a number of, at most,
$C|||z|||^2log(|||z|||)$ different arcs $\frac12\gamma$ needed to
cover $\partial\Omega^2$. Then an iteration of the inequality above
yields
\begin{eqnarray*}
e^{-C|||z|||^2
log(|||z|||)}\leq|\frac{\phi'(\tau_1)}{\phi'(\tau_2)}| \leq
e^{C|||z|||^2 log(|||z|||)}
\end{eqnarray*}
for any two arbitrary points $\tau_1,\tau_2\in\partial\Omega^2$,
allowing us to finish the proof in the case $\mathcal{H}^2$. The
transformation $z \rightarrow 1/(z - z_0)$ where, as before,
$z_0\in\Omega^2$, $dist(z_0,\partial\Omega^2)>> 1/|||z|||$, allows
us to reduce the estimate for $\mathcal{H}^1$ to the previous case.

%%%%%%%%%%%%%%%%%%%%%%%%%%%%%%%%%%%%%%%%%%%%%%%%%%%%%%%%%%%%%%%%%%%%%%%%%%%%%
%%%%%%%%%%%%%%%%%%%%%%%%%%%%%%%%%%%%%%%%%%%%%%%%%%%%%%%%%%%%%%%%%%%%%%%%%%%%%

\section{Preliminary estimates}

The following subsection are devoted to show the regularity of the
different elements involved in the problem: the Birkhoff-Rott
integral, $z_t(\al,t)$, $\varpi_t(\al,t)$, $\varpi(\al,t)$; the
difference of the gradient of the pressure in the normal direction
$\sigma(\al,t)$ and its time derivative $\sigma_t(\al,t)$. We shall
concentrate our attention in  the case of a closed contour, because
for a periodic domain in the horizontal space variable the treatment
is completely analogous (see \cite{DY2}).

\subsection{Estimates for $BR(z,\varpi)$}

%%%%%%%%%%%%%%%%%%%%%%%%%%%%%%%%%%%%%%%%%%%%%%%%%%%%%%%%%%%%%%%%%%%%%%%%%%%%%
%%%%%%%%%%%%%%%%%%%%%%%%%%%%%%%%%%%%%%%%%%%%%%%%%%%%%%%%%%%%%%%%%%%%%%%%%%%%%

In this section we show that the Birkhoff-Rott integral is
as regular as $\da z$.

\begin{lemma}The following estimate holds
\begin{eqnarray}\label{nsibr}
\|BR(z,\varpi)\|_{H^k}\leq
C(\|\F(z)\|^2_{L^\infty}+\|z\|^2_{H^{k+1}}+\|\varpi\|^2_{H^{k}})^j,
\end{eqnarray}
for $k\geq 2$, where $C$ and $j$ are constants independent of $z$
and $\varpi$.
\end{lemma}
\begin{rem}
Using this estimate for $k=2$ we find easily that
\begin{eqnarray}\label{nliibr}
\|\da BR(z,\varpi)\|_{L^\infty}\leq
C(\|\F(z)\|^2_{L^\infty}+\|z\|^2_{H^{3}}+\|\varpi\|^2_{H^{2}})^j,
\end{eqnarray}
which shall be used through out the paper.
\end{rem}

Proof: We shall present the proof for $k=2$. Let us write
\begin{align*}
\begin{split}
BR(z,\varpi)(\al,t)&=\frac{1}{2\pi}\int_{-\pi}^\pi C_1(\al,\beta)\varpi(\al-\beta)d\beta
+\frac{\dpa z(\al)}{2|\da z(\al)|^2}H(\varpi)(\al)
\end{split}
\end{align*}
where $C_1$ is given by \eqref{fb}. The boundedness of the term $C_1$ in
$L^\infty$ gives us easily
\begin{equation}\label{brl2}
\|BR(z,\varpi)\|_{L^2}\leq C\|\F(z)\|^2_{L^\infty}\|z\|^2_{C^2}\|\varpi\|_{L^2}.
\end{equation}
 In $\da^2BR(z,\varpi)$, the most singular terms are given by
$$
P_1(\al)=\frac{1}{2\pi}PV\int_{-\pi}^\pi\da^2\varpi(\al-\beta)\frac{(z(\al)-z(\al-\beta))^{\bot}}{|z(\al)-z(\al-\beta)|^2}d\beta,
$$
$$
P_2(\al)=\frac{1}{2\pi}PV\int_{-\pi}^\pi\varpi(\al-\beta)\frac{\da^2z(\al)-\da^2z(\al-\beta)}{|z(\al)-z(\al-\beta)|^2}d\beta,
$$
$$
P_3(\al)=-\frac{1}{\pi}PV\int_{-\pi}^\pi\varpi(\al-\beta)\frac{(z(\al)-z(\al-\beta))^\bot}{|z(\al)-z(\al-\beta)|^4}\big(z(\al)-z(\al-\beta))\cdot(\da^2 z(\al)-\da^2 z(\al-\beta))\big)d\beta.
$$
Again we have the expression
\begin{align*}
\begin{split}
P_1(\al)&=\frac{1}{2\pi}\int_{-\pi}^\pi C_1(\al,\beta)\da^2\varpi(\al-\beta)d\beta
+\frac{\dpa
z(\al)}{2|\da z(\al)|^2}H(\da^2\varpi)(\al)d\al,
\end{split}
\end{align*} giving us
\begin{align}\label{q1nl22dbr}
\begin{split}
|P_1(\al)|&\leq C\|\F(z)\|^j_{L^\infty}\|z\|^j_{C^2}(\|\da^2\varpi\|_{L^2}+|H(\da^2\varpi)(\alpha)|).
\end{split}
\end{align}
Next let us write $P_2=Q_1+Q_2+Q_3$ where
$$
Q_1(\al)=\frac{1}{2\pi}\int_{-\pi}^\pi(\varpi(\al-\beta)-\varpi(\al))\frac{\da^2 z(\al)-\da^2 z(\al-\beta)}{|z(\al)-z(\al-\beta)|^2}d\beta,
$$
$$
Q_2(\al)=\frac{\varpi(\al)}{2\pi}\int_{-\pi}^\pi(\da^2z(\al)-\da^2z(\al-\beta))\big(\frac{1}{|z(\al)-z(\al-\beta)|^2}-\frac{1}{|\da z(\al)|^2|\beta|^2}\big)d\beta,
$$
$$
Q_3(\al)=\frac{1}{2\pi}\frac{\varpi(\al)}{|\da z(\al)|^2}\int_{-\pi}^\pi\!\!(\da^2z(\al)\!-\!\da^2z(\al\!-\!\beta))\big(\frac{1}{|\beta|^2}\!-\!\frac{1}{4\sin^2(\beta/2)}\big)d\beta\!+\!\frac{1}{2}\frac{\varpi(\al)}{|\da z(\al)|^2}\la (\da^2z)(\al),
$$
where $\la=\da H$.

Using that
$$|\da^2z(\al)-\da^2z(\al-\beta)|\leq|\beta|^\delta\|z\|_{C^{2,\delta}},$$
we get $|Q_1(\al)|+|Q_2(\al)|\leq
\|\varpi\|_{C^1}\|\F(z)\|^{j}\|z\|^j_{C^{2,\delta}},$ while for
$Q_3$ we have
\begin{align*}
|Q_3(\al)|&\leq C\|\varpi\|_{L^\infty}\|\F(z)\|_{L^\infty}(\|z\|_{C^2}+|\la (\da^2z)(\al)|),
\end{align*}
that is
\begin{equation}\label{q2nl22dbr}
|P_2(\al)|\leq (1+|\la (\da^2z)(\al)|)\|\varpi\|_{C^1}\|\F(z)\|^{j}\|z\|^j_{C^{2,\delta}}.
\end{equation}
Let us now consider  $P_3=Q_4+Q_5+Q_6+Q_7+Q_8+Q_9$, where
\begin{align*}
Q_4=\frac{-1}{\pi}\int_{-\pi}^\pi(\varpi(\al\!-\!\beta)\!-\!\varpi(\al))\frac{(z(\al)\!-\!z(\al\!-\!\beta))^\bot}{|z(\al)\!-\!z(\al\!-\!\beta)|^4}\big((z(\al)\!-\!z(\al\!-\!\beta))\!\cdot\!(\da^2 z(\al)\!-\!\da^2 z(\al\!-\!\beta))\big)d\beta,
\end{align*}
$$Q_5=-\frac{\varpi(\al)}{\pi}\int_{-\pi}^\pi\frac{(z(\al)\!-\!z(\al\!-\!\beta)\!-\!\da z(\al)\beta)^{\bot}}{|z(\al)\!-\!z(\al\!-\!\beta)|^4}\big((z(\al)\!-\!z(\al\!-\!\beta))\cdot(\da^2 z(\al)\!-\!\da^2 z(\al\!-\!\beta))\big)d\beta,
$$
$$
Q_6=-\frac{\varpi(\al)\dpa z(\al)}{\pi}\int_{-\pi}^\pi\frac{\beta(z(\al)\!-\!z(\al\!-\!\beta)\!-\!\da z(\al)\beta)\cdot(\da^2 z(\al)\!-\!\da^2 z(\al\!-\!\beta))}{|z(\al)\!-\!z(\al\!-\!\beta)|^4}d\beta,
$$
$$
Q_7=-\frac{\varpi(\al)\dpa z(\al)}{\pi}\da z(\al)\cdot\!\!\int_{-\pi}^\pi\beta^2(\da^2 z(\al)-\da^2 z(\al-\beta))\big(\frac{1}{|z(\al)\!-\!z(\al\!-\!\beta)|^4}-\frac{1}{|\da z(\al)|^4|\beta|^4}\big)d\beta,
$$
$$
Q_8=-\frac{\varpi(\al)\dpa z(\al)}{\pi|\da z(\al)|^4}\da z(\al)\cdot\!\!\int_{-\pi}^\pi(\da^2 z(\al)-\da^2 z(\al-\beta))\big(\frac{1}{|\beta|^2}-\frac{1}{4\sin^2(\beta/2)}\big)d\beta,
$$
and
$$
Q_9=-\frac{\varpi(\al)\dpa z(\al)}{|\da z(\al)|^4}\da z(\al)\cdot \la(\da^2 z(\al)).
$$
Proceeding as before we get
\begin{equation*}
|P_3(\al)|\leq C(1+|\la (\da^2z)(\al)|)\|\varpi\|_{C^1}\|\F(z)\|^j_{L^\infty}\|z\|^j_{C^{2,\delta}},
\end{equation*}
which together with \eqref{q1nl22dbr} and \eqref{q2nl22dbr} gives
us the estimate
$$
|(P_1+P_2+P_3)(\al)|\leq C(1+|\la (\da^2z)(\al)|+|H(\da^2\varpi)(\alpha)|)\|\varpi\|_{C^1}(\|\F(z)\|^j_{L^\infty}+\|z\|^j_{H^3}).
$$
For the rest of the terms in $\da^2 BR(z,\varpi)$ we obtain
analogous estimates allowing us to conclude the equality
\begin{equation*}\label{enl22dbr}
\|\da^2 BR(z,\varpi)\|_{L^2}\leq C(1+\|\da^3z\|_{L^2}+\|\da^2\varpi\|_{L^2}) \|\varpi\|_{C^1}\|\F(z)\|^j_{L^\infty}\|z\|^j_{C^{2,\delta}}.
\end{equation*}
Finally the Sobolev inequalities yield \eqref{nsibr} for $k=2$.

%%%%%%%%%%%%%%%%%%%%%%%%%%%%%%%%%%%%%%%%%%%%%%%%%%%%%%%%%%%%%%%%%%%%%%%%%%%%%
%%%%%%%%%%%%%%%%%%%%%%%%%%%%%%%%%%%%%%%%%%%%%%%%%%%%%%%%%%%%%%%%%%%%%%%%%%%%%

\subsection{Estimates for $z_t(\al,t)$}

%%%%%%%%%%%%%%%%%%%%%%%%%%%%%%%%%%%%%%%%%%%%%%%%%%%%%%%%%%%%%%%%%%%%%%%%%%%%%
%%%%%%%%%%%%%%%%%%%%%%%%%%%%%%%%%%%%%%%%%%%%%%%%%%%%%%%%%%%%%%%%%%%%%%%%%%%%%

This section is devoted to show that $z_t$ is as regular as $\da z$.

\begin{lemma}The following estimate holds
\begin{eqnarray}\label{nszt}
\|z_t\|_{H^k}\leq
C(\|\F(z)\|^2_{L^\infty}+\|z\|^2_{H^{k+1}}+\|\varpi\|^2_{H^{k}})^j,
\end{eqnarray}
for $k\geq 2$.
\end{lemma}

Proof: It follows easily from formulas  \eqref{eev}, \eqref{fc2}
together with the estimates obtained in the last section.

%%%%%%%%%%%%%%%%%%%%%%%%%%%%%%%%%%%%%%%%%%%%%%%%%%%%%%%%%%%%%%%%%%%%%%%%%%%%%
%%%%%%%%%%%%%%%%%%%%%%%%%%%%%%%%%%%%%%%%%%%%%%%%%%%%%%%%%%%%%%%%%%%%%%%%%%%%%

\subsection{Estimates for $\varpi_t$}

%%%%%%%%%%%%%%%%%%%%%%%%%%%%%%%%%%%%%%%%%%%%%%%%%%%%%%%%%%%%%%%%%%%%%%%%%%%%%
%%%%%%%%%%%%%%%%%%%%%%%%%%%%%%%%%%%%%%%%%%%%%%%%%%%%%%%%%%%%%%%%%%%%%%%%%%%%%

This section is devoted to show that $\varpi_t$ is as regular as
$\da \varpi$

\begin{lemma}The following estimate holds
\begin{eqnarray}\label{nswt}
\|\varpi_t\|_{H^k}\leq C
\exp(C|||z|||^p)(\|\F(z)\|^2_{L^\infty}+\|z\|^2_{H^{k+2}}+\|\varpi\|^2_{H^{k+1}}+\|\varphi\|^2_{H^{k+1}})^j,
\end{eqnarray}
for $k\geq 1$.
\end{lemma}

Proof: In the following we shall work the details of the proof only
when  $k=1$, since the cases $k\geq 2$ can be treated analogously.
Formula \eqref{nose} yields
\begin{align}
\begin{split}\label{nose2}
\varpi_t(\al,t)+T(\varpi_t)(\al,t)&=I_1(\al,t)+I_2(\al,t)
-2\varphi(\al,t)\da\varphi(\al,t)+R(\al,t),
\end{split}
\end{align}
where
$$
I_1=\frac{-1}{\pi}\int_{-\pi}^\pi\frac{(z_t(\al)-z_t(\al-\beta))^{\bot}\cdot\da z(\al)}{|z(\al)-z(\al-\beta)|^2}\varpi(\al-\beta)d\beta,
$$
$$
I_2=\frac{2}{\pi}\int_{-\pi}^\pi\frac{(z(\al)\!-\!z(\al\!-\!\beta))^{\bot}\cdot\da z(\al)}{|z(\al)\!-\!z(\al\!-\!\beta)|^4}(z(\al)\!-\!z(\al\!-\!\beta))\cdot(z_t(\al)\!-\!z_t(\al\!-\!\beta))\varpi(\al\!-\!\beta)d\beta,
$$
and
$$
R=\frac{c(\al,t)}{\pi}\int_{-\pi}^\pi \!\da z(\al,t)\cdot \da
BR(z,\varpi)(\al,t) d\al+2\mathrm{g}\da z_2(\al,t).
$$

>From Theorem \ref{oil2l2} we get
$$
\|\varpi_t\|_{L^2}\leq \|(I+T)^{-1}\|_{L^2\rightarrow L^2}(\|I_1\|_{L^2}+\|I_1\|_{L^2}+2\|\varphi\da\varphi\|_{L^2}+\|R\|_{L^2}),
$$
and proceeding as before, using the estimates above, we obtain
\begin{equation}\label{nl2wt}
\|w_t\|_{L^2}\leq\exp(C|||z|||^p)(\|\F(z)\|^2_{L^\infty}+\|z\|^2_{H^{3}}+\|\varpi\|^2_{H^{2}}+\|\varphi\|^2_{H^{2}})^j.
\end{equation}
Next we shall show that in the singular case we have:
\begin{equation}\label{nl2dwt}
\|\da w_t\|_{L^2}\leq\exp(C|||z|||^p)(\|\F(z)\|^2_{L^\infty}+\|z\|^2_{H^{3}}+
\|\varpi\|^2_{H^{2}}+\|\varphi\|^2_{H^{2}})^j.
\end{equation}
To see it let us  take a derivative in \eqref{nose2} to obtain the
identity
\begin{align}
\begin{split}\label{nose3}
\da\varpi_t(\al,t)+T(\da\varpi_t)(\al,t)&=J_1(\al,t)+J_2(\al,t)+J_3(\al,t)+\da I_1(\al,t)+\da I_2(\al,t)\\
&\quad-\da^2(\varphi^2)(\al,t)+\da R(\al,t),
\end{split}
\end{align}
where
$$
J_1=\frac{-1}{\pi}\int_{-\pi}^\pi\frac{(\da z(\al)-\da z(\al-\beta))^{\bot}\cdot\da z(\al)}{|z(\al)-z(\al-\beta)|^2}\varpi_t(\al-\beta)d\beta,
$$
$$
J_2=\frac{2}{\pi}\int_{-\pi}^\pi\frac{(z(\al)\!-\!z(\al\!-\!\beta))^{\bot}\!\cdot\!\da z(\al)}{|z(\al)\!-\!z(\al\!-\!\beta)|^4}(z(\al)\!-\!z(\al\!-\!\beta))\!\cdot\!(\da z(\al)\!-\!\da z(\al\!-\!\beta))\varpi_t(\al\!-\!\beta)d\beta,
$$
and
$$
J_3=\frac{-1}{\pi}\int_{-\pi}^\pi\frac{(z(\al)-z(\al-\beta))^{\bot}\cdot\da^2 z(\al)}{|z(\al)-z(\al-\beta)|^2}\varpi_t(\al-\beta)d\beta,
$$

Using Theorem \ref{oil2l2} in \eqref{nose3} we get
$$
\|\da\varpi_t\|_{L^2}\leq \|(I+T)^{-1}\|_{L^2\rightarrow L^2}(\sum_{l=1}^3\|J_l\|_{L^2}+\|\da I_1\|_{L^2}+\|\da I_2\|_{L^2}+\|\da^2(\varphi^2)\|_{L^2}+\|\da R\|_{L^2}),
$$
A straightforward calculation yields
$$
\|\da^2(\varphi^2)\|_{L^2}+\|\da R\|_{L^2}\leq C(\|\F(z)\|^2_{L^\infty}+\|z\|^2_{H^{3}}+
\|\varpi\|^2_{H^{2}}+\|\varphi\|^2_{H^{2}})^j.
$$
To estimate the other terms we write:
$$
J_1=\frac{-1}{\pi}\int_{-\pi}^\pi C_2(\al,\beta)\varpi_t(\al-\beta)d\beta-\frac{(\da^2z(\al))^{\bot}\cdot\da z(\al)}{|\da z(\al)|^2}H(\varpi_t)(\al),
$$
where
$$
C_2(\al,\beta)=[\frac{(\da z(\al)-\da z(\al-\beta))^{\bot}\cdot\da z(\al)}{|z(\al)-z(\al-\beta)|^2}-\frac{(\da^2z(\al))^{\bot}\cdot\da z(\al)}{|\da z(\al)|^22\tan(\beta/2)}].
$$
Then
$$
|J_1(\al)|\leq C\|\F(z)\|^k_{L^\infty}\|z\|^k_{C^{2,\delta}}(\int_{-\pi}^\pi |\beta|^{\delta-1}|\varpi_t(\al-\beta)|d\beta+|H(\varpi_t)(\al)|),
$$
and using \eqref{nl2wt} we have
$$
\|J_1\|_{L^2}\leq \exp(C|||z|||^p)(\|\F(z)\|^2_{L^\infty}+\|z\|^2_{H^{3}}+
\|\varpi\|^2_{H^{2}}+\|\varphi\|^2_{H^{2}})^j.
$$
Next we  rewrite $J_2$ as follows
$$
\frac{2}{\pi}\int_{-\pi}^\pi\frac{(z(\al)\!-\!z(\al\!-\!\beta)\!-\!\da z(\al)\beta)^{\bot}\!\cdot\!\da z(\al)}{|z(\al)\!-\!z(\al\!-\!\beta)|^4}(z(\al)\!-\!z(\al\!-\!\beta))\!\cdot\!(\da z(\al)\!-\!\da z(\al\!-\!\beta))\varpi_t(\al\!-\!\beta)d\beta,
$$
which is a more regular term than $J_1$. Since $J_3$ is also more
regular than $J_1$ we finally  get
$$
\|J_2\|_{L^2}+\|J_3\|_{L^2}\leq \exp(C|||z|||^p)(\|\F(z)\|^2_{L^\infty}+\|z\|^2_{H^{3}}+
\|\varpi\|^2_{H^{2}}+\|\varphi\|^2_{H^{2}})^j.
$$
The most singular term in $\da I_1$ is given by
$$
K_1=\frac{-1}{\pi}\int_{-\pi}^\pi\frac{(\da z_t(\al)-\da z_t(\al-\beta))^{\bot}\cdot\da z(\al)}{|z(\al)-z(\al-\beta)|^2}\varpi(\al-\beta)d\beta,
$$
 and will be estimated using the following splitting $K_1=L_1+L_2+L_3+L_4$
where
$$
L_1=\frac{1}{\pi}\int_{-\pi}^\pi\frac{(\da z_t(\al)-\da z_t(\al-\beta))^{\bot}\cdot\da z(\al)}{|z(\al)-z(\al-\beta)|^2}
(\varpi(\al)-\varpi(\al-\beta))d\beta,
$$
$$
L_2=\frac{-\varpi(\al)}{\pi}\int_{-\pi}^\pi(\da z_t(\al)-\da z_t(\al-\beta))^{\bot}\cdot\da z(\al)[\frac{1}{|z(\al)-z(\al-\beta)|^2}-\frac1{|\da z(\al)|^2|\beta|^2}]d\beta,
$$
$$
L_3=\frac{-1}{\pi}\frac{\varpi(\al)\da z(\al)}{|\da z(\al)|^2}\cdot\int_{-\pi}^\pi(\da z_t(\al)-\da z_t(\al-\beta))^{\bot}[\frac1{\beta^2}-\frac1{4\sin^2(\beta/2)}]d\beta,
$$
and
$$
L_4=\frac{-1}{\pi}\frac{\varpi(\al)\da z(\al)}{|\da z(\al)|^2}\cdot\la(\da z_t).
$$
Since $|\da z_t(\al)-\da z_t(\al-\beta)|\leq
|\beta|\int_0^1|\da^2z_t(\al+(s-1)\beta)|ds$ we have
$$
|K_1|\leq C\|\F(z)\|^k_{L^\infty}\|z\|^k_{H^{3}}\|\varpi\|_{C^1}(\int_0^1|\da^2z_t(\al+(s-1)\beta)|ds+|\la(\da z_t)(\al)|).
$$ From \eqref{nszt} we obtain the estimates
$$
\|K_1\|_{L^2}\leq C(\|\F(z)\|^2_{L^\infty}+\|z\|^2_{H^{3}}+
\|\varpi\|^2_{H^{2}})^j
$$
and
$$
\|\da I_1\|_{L^2}\leq C(\|\F(z)\|^2_{L^\infty}+\|z\|^2_{H^{3}}+
\|\varpi\|^2_{H^{2}})^j.
$$
Next we rewrite $I_2$ in the form
$$
\frac{2}{\pi}\int_{-\pi}^\pi\frac{(z(\al)\!-\!z(\al\!-\!\beta)-\da z(\al)\beta)^{\bot}\!\cdot\!\da z(\al)}{|z(\al)\!-\!z(\al\!-\!\beta)|^4}(z(\al)\!-\!z(\al\!-\!\beta))\!\cdot\!(z_t(\al)\!-\!z_t(\al\!-\!\beta))\varpi(\al\!-\!\beta)d\beta,
$$
which shows that $I_2$ is more regular than $I_1$ and, therefore,
the estimate for  $\da I_2$ follow easily with the same methods that
we used with $\da I_1$, allowing us to finish the proof.

%%%%%%%%%%%%%%%%%%%%%%%%%%%%%%%%%%%%%%%%%%%%%%%%%%%%%%%%%%%%%%%%%%%%%%%%%%%%%
%%%%%%%%%%%%%%%%%%%%%%%%%%%%%%%%%%%%%%%%%%%%%%%%%%%%%%%%%%%%%%%%%%%%%%%%%%%%%

\subsection{Estimates for $\varpi$}

%%%%%%%%%%%%%%%%%%%%%%%%%%%%%%%%%%%%%%%%%%%%%%%%%%%%%%%%%%%%%%%%%%%%%%%%%%%%%
%%%%%%%%%%%%%%%%%%%%%%%%%%%%%%%%%%%%%%%%%%%%%%%%%%%%%%%%%%%%%%%%%%%%%%%%%%%%%

In this section we show that the amplitude of the vorticity $\varpi$
lies at the same level than $\da z$. We shall consider $z\in
H^k(\T)$, $\varphi\in H^{k-\frac12}(\T)$ and $\varpi\in H^{k-2}(\T)$
as part of the energy estimates. The inequality below yields
$\varpi\in H^{k-1}(\T)$.

\begin{lemma}The following estimate holds
\begin{eqnarray}\label{nsw}
\|\varpi \|_{H^k}\leq
C(\|\F(z)\|^2_{L^\infty}+\|z\|^2_{H^{k+1}}+\|\varpi\|^2_{H^{k-1}}+\|\varphi\|^2_{H^{k}})^j,
\end{eqnarray}
for $k\geq 2$.
\end{lemma}

Proof: We shall present  the proof for $k=2$, being the rest of the
cases completely analogous. Since $\varpi=2|\da z|\varphi+2|\da
z|^2c$  the identity $|\da z|^2=A(t)$ gives us the equality
\begin{align*}
\da^2\varpi(\al)&=2|\da z(\al)|\da^2\varphi(\al)-\da(2\da z\cdot\da
BR(z,\varpi))(\al),
\end{align*}
from which we easily get
$$
\|\da^2\varpi\|_{L^2}\leq
2\|z\|_{C^1}\|\da^2\varphi\|_{L^2}+\|\da(2\da z\cdot\da
BR(z,\varpi))\|_{L^2}.
$$
Therefore in order to get the estimate \eqref{nsw} for $k=2$ we need
to show that the following inequality holds
\begin{equation}\label{hec}
\|\da(2\da z\cdot\da BR(z,\varpi))\|_{L^2}\leq
C\|\F(z)\|^j_{L^{\infty}}\|z\|^j_{H^3}\|\varpi\|_{H^1}.
\end{equation}

To see that we can write
$$
2\da z(\al)\cdot\da BR(z,\varpi)(\al)=T(\da\varpi)(\al)+R_1(\al)+R_2(\al),
$$
where
$$
R_1(\al)=\frac{1}{\pi}\da z(\al)\cdot\int_{-\pi}^\pi\frac{(\da z(\al)-\da z(\al-\beta))^{\bot}}{|z(\al)-z(\al-\beta)|^2}\varpi(\al-\beta)d\beta,
$$
and
$$
R_2(\al)=-\frac{2}{\pi}\da z(\al)\cdot\!\!\int_{-\pi}^\pi\frac{(z(\al)\!-\! z(\al\!-\!\beta))^{\bot}}{|z(\al)\!-\!z(\al\!-\!\beta)|^4}(z(\al)\!-\!z(\al\!-\!\beta))\cdot(\da z(\al)\!-\!\da z(\al\!-\!\beta))\varpi(\al\!-\!\beta).
$$
Then we have $\|T(\da\varpi)\|_{H^1}\leq
C\|\F(z)\|^4_{L^{\infty}}\|z\|^4_{C^{2,\delta}}\|\da \varpi\|_{L^2}$
from \eqref{ol2h1}, so that we only need to estimate $\da R_1$ and
$\da R_2$ in $L^2$ to get \eqref{hec}.

Next we consider the most singular terms in $\da R_1$, namely:
$$
S_1(\al)=\frac{1}{\pi}\da z(\al)\cdot\int_{-\pi}^\pi\frac{(\da^2 z(\al)-\da z^2(\al-\beta))^{\bot}}{|z(\al)-z(\al-\beta)|^2}\varpi(\al-\beta)d\beta,
$$
$$
S_2(\al)=\frac{1}{\pi}\da z(\al)\cdot\int_{-\pi}^\pi\frac{(\da
z(\al)-\da
z(\al-\beta))^{\bot}}{|z(\al)-z(\al-\beta)|^2}\da\varpi(\al-\beta)d\beta,
$$
and we use the decomposition
\begin{align*}
S_2(\al)&=\frac{1}{\pi}\da z(\al)\cdot\int_{-\pi}^\pi[\frac{(\da z(\al)-\da z(\al-\beta))^{\bot}}{|z(\al)-z(\al-\beta)|^2}-\frac{(\da^2z(\al))^{\bot}}{|\da z(\al)|^22\tan(\beta/2)}]\da\varpi(\al-\beta)d\beta\\
&\quad-\frac{\dpa z(\al)\cdot\da^2z(\al)}{|\da z(\al)|^2}H(\da\varpi)(\al).
\end{align*}
to obtain
$$
|S_2(\al)|\leq
C(\|\F(z)\|^j_{L^{\infty}}\|z\|^j_{C^{2,\delta}}(\|\da
\varpi\|_{L^2}+|H(\da\varpi)(\al)|+\int_{-\pi}^\pi|\beta|^{\delta-1}|\da\varpi(\al-\beta)|d\beta),
$$
that is $\|S_2\|_{L^2}\leq
C\|\F(z)\|^j_{L^{\infty}}\|z\|^j_{C^{2,\delta}}\|\da
\varpi\|_{L^2}.$

In $S_1$ we have the splitting $U_1+U_2+U_3+U_4$ where
$$
U_1(\al)=\frac{1}{\pi}\da z(\al)\cdot\int_{-\pi}^\pi\frac{(\da^2 z(\al)-\da z^2(\al-\beta))^{\bot}}{|z(\al)-z(\al-\beta)|^2}(\varpi(\al-\beta)-\varpi(\al))d\beta,
$$
$$
U_2(\al)=\frac{1}{\pi}\varpi(\al)\da z(\al)\cdot\int_{-\pi}^\pi(\da^2 z(\al)\!-\!\da z^2(\al\!-\!\beta))^{\bot}[\frac{1}{|z(\al)\!-\!z(\al\!-\!\beta)|^2}\!-\!\frac{1}{|\da z(\al)|^2\beta^2}]d\beta,
$$
$$
U_3(\al)=\frac{1}{\pi}\varpi(\al)\frac{\da z(\al)}{|\da z(\al)|^2}\cdot\int_{-\pi}^\pi(\da^2 z(\al)-\da z^2(\al-\beta))^{\bot}[\frac{1}{\beta^2}-\frac{1}{4\sin^2(\beta/2)}]d\beta,
$$
and
$$
U_4(\al)=\varpi(\al)\frac{\da z(\al)}{|\da z(\al)|^2}\cdot\la(\da^2 z)(\al).
$$
Then in $U_1$ we use the identity
\begin{equation}\label{hlp}
\da^2 z(\al)-\da z^2(\al-\beta)=\beta\int_0^1\da^3z(\al+(s-1)\beta)ds
\end{equation}
to get
$$
|U_1(\al)|\leq C\|\F(z)\|^2_{L^{\infty}}\|z\|_{C^{2,\delta}}\|\varpi\|_{C^\delta}\int_0^1\int_{-\pi}^\pi|\beta|^{\delta-1}|\da^3z(\al+(s-1)\beta)|d\beta ds,
$$
and therefore $\|U_1\|_{L^2}\leq C\|\F(z)\|^2_{L^{\infty}}\|z\|^2_{H^3}\|\varpi\|_{H^1}$.

To estimate $U_2$ and $U_3$ we can use again \eqref{hlp}. For $U_4$
the control is easier.

To finish the argument we rewrite $R_2$ as follows:
$$
-\frac{2}{\pi}\da z(\al)\cdot\!\!\int_{-\pi}^\pi\frac{(z(\al)\!-\!
z(\al\!-\!\beta)\!-\!\da
z(\al)\beta)^{\bot}}{|z(\al)\!-\!z(\al\!-\!\beta)|^4}(z(\al)\!-\!z(\al\!-\!\beta))\cdot(\da
z(\al)\!-\!\da z(\al\!-\!\beta))\varpi(\al\!-\!\beta)d\beta,
$$
expressing the fact that with the same method, $\da R_2$ is easier
to estimate than $\da R_1$.

%%%%%%%%%%%%%%%%%%%%%%%%%%%%%%%%%%%%%%%%%%%%%%%%%%%%%%%%%%%%%%%%%%%%%%%%%%%%%
%%%%%%%%%%%%%%%%%%%%%%%%%%%%%%%%%%%%%%%%%%%%%%%%%%%%%%%%%%%%%%%%%%%%%%%%%%%%%

\subsection{Estimates for $\sigma$}

%%%%%%%%%%%%%%%%%%%%%%%%%%%%%%%%%%%%%%%%%%%%%%%%%%%%%%%%%%%%%%%%%%%%%%%%%%%%%
%%%%%%%%%%%%%%%%%%%%%%%%%%%%%%%%%%%%%%%%%%%%%%%%%%%%%%%%%%%%%%%%%%%%%%%%%%%%%

Here we prove that   $\sigma$, the difference of the gradient of the
pressure in the normal direction,  is at the same level than $\da^2
z$.

\begin{lemma}The following estimate holds
\begin{eqnarray}\label{nswbis}
\|\sigma \|_{H^k}\leq
C\exp(C|||z|||^p)(\|\F(z)\|^2_{L^\infty}+\|z\|^2_{H^{k+2}}+\|\varpi\|^2_{H^{k+1}}+\|\varphi\|^2_{H^{k+1}})^j,
\end{eqnarray}
for $k\geq 2$.
\end{lemma}
Proof: We shall give the details of the case  $k=2$. Let us recall
the formula for $\sigma(\al)$:
\begin{equation}\label{fs2}
\frac{\sigma}{\rho_2}=(\partial_t BR(z,\varpi)\!+\!\frac{\varphi}{|\da z|}\da BR(z,\varpi))\cdot \dpa z\!+\!\frac12\frac{\varpi}{|\da z|^2}(\da
z_t\!+\!\frac{\varphi}{|\da z|}\da^2 z)\cdot
\dpa z\!+\!\mathrm{g}\da z_1.
\end{equation}
then from previous sections we have:
$$
\|\sigma \|_{L^2}\leq
C\exp(C|||z|||^p)(\|\F(z)\|^2_{L^\infty}+\|z\|^2_{H^{4}}+\|\varpi\|^2_{H^{3}}+\|\varphi\|^2_{H^{3}})^j.
$$
To control $\|\da^2\sigma \|_{L^2}$ we only have to deal with
$\da^2(\partial_t BR(z,\varpi)\cdot \dpa z)$, because the remainder
terms have been already estimated. Again we shall consider the most
singular parts:
$$
I_1=\frac{1}{2\pi}\int_{-\pi}^\pi \frac{(z(\al)\!-\!z(\al\!-\!\beta))\cdot\da
z(\al)}{|z(\al)-z(\al-\beta)|^2}\da^2\varpi_t(\al-\beta)d\beta,
$$
$$
I_2=\frac{1}{2\pi}\int_{-\pi}^\pi \frac{(\da^2z_t(\al)\!-\!\da^2z_t(\al\!-\!\beta))\cdot\da
z(\al)}{|z(\al)-z(\al-\beta)|^2}\varpi(\al-\beta)d\beta,
$$
$$
I_3=\frac{-1}{\pi}\int_{-\pi}^\pi \frac{(z(\al)\!-\!z(\al\!-\!\beta))\!\cdot\!\da
z(\al)}{|z(\al)\!-\!z(\al\!-\!\beta)|^4}(z(\al)\!-\!z(\al\!-\!\beta))\!\cdot\!(\da^2z_t(\al)\!-\!\da^2z_t(\al\!-\!\beta))\varpi(\al-\beta)d\beta.
$$
We have
$$
I_1=\frac{1}{2\pi}\int_{-\pi}^\pi E(\al,\beta)\da^2\varpi_t(\al-\beta)d\beta+\frac12H(\da^2\varpi_t)(\al),
$$
where
$$
E(\al,\beta)=\frac{(z(\al)\!-\!z(\al\!-\!\beta))\cdot\da
z(\al)}{|z(\al)-z(\al-\beta)|^2}-\frac{1}{2\tan(\beta/2)}.
$$
Since $\|E\|_{L^\infty}\leq C\|\F(z)\|^2_{L^{\infty}}\|z\|^2_{C^2}$
we can estimate $I_1$ throughout inequality \eqref{nswt}.

The  equality
\begin{equation*}
\da^2 z_t(\al)-\da^2 z_t(\al-\beta)=\beta\int_0^1\da^3z_t(\al+(s-1)\beta)ds
\end{equation*}
let us to get
$$
|I_2|+|I_3|\leq C\|\F(z)\|^2_{L^{\infty}}\|z\|^2_{C^2}\|\varpi\|_{C^1}(\int_0^1\int_{-\pi}^\pi|\da^3z_t(\al+(s-1)\beta)|ds+|\la(\da^2z_t)(\al)|)
$$
and \eqref{nszt} take care of the rest.

%%%%%%%%%%%%%%%%%%%%%%%%%%%%%%%%%%%%%%%%%%%%%%%%%%%%%%%%%%%%%%%%%%%%%%%%%%%%%
%%%%%%%%%%%%%%%%%%%%%%%%%%%%%%%%%%%%%%%%%%%%%%%%%%%%%%%%%%%%%%%%%%%%%%%%%%%%%

\subsection{Estimate for $\sigma_t$}

%%%%%%%%%%%%%%%%%%%%%%%%%%%%%%%%%%%%%%%%%%%%%%%%%%%%%%%%%%%%%%%%%%%%%%%%%%%%%
%%%%%%%%%%%%%%%%%%%%%%%%%%%%%%%%%%%%%%%%%%%%%%%%%%%%%%%%%%%%%%%%%%%%%%%%%%%%%

In this section we obtain an upper bound for the $L^\infty$ norm of
$\sigma_t$ that will be used in the energy inequalities and  in the
treatment of the Rayleigh-Taylor condition.

\begin{lemma} The following estimate holds
\begin{eqnarray}\label{nlinftyst}
\|\sigma_t \|_{L^\infty}\leq
C\exp(C|||z|||^p)(\|\F(z)\|^2_{L^\infty}+\|z\|^2_{H^4}+\|\varpi\|^2_{H^3}+\|\varphi\|^2_{H^3})^j.
\end{eqnarray}
\end{lemma}

Proof: Let us consider \eqref{fs2} the splitting
$\sigma/\rho_2=P_1+P_2+P_3+P_4+P_5$ where
$$
P_1=\partial_t BR(z,\varpi)\cdot\dpa z,\quad P_2=\frac{\varphi}{|\da z|}\da BR(z,\varpi)\cdot\dpa z,$$
$$
P_3=\frac12\frac{\varpi}{|\da z|^2}\da
z_t\cdot\dpa z,\quad P_4=\frac{\varphi}{|\da z|}\da^2 z\cdot\dpa z,\quad P_5=\mathrm{g}\da z_1.
$$
Estimate \eqref{nszt} yields $\|\dpt
P_5\|_{L^{\infty}}\leq \|\mathrm{g}\dpt\da z_1\|_{H^1}\leq
C(\|\F(z)\|^2_{L^\infty}+\|z\|^2_{H^3}+\|\varpi\|^2_{H^2})^j.$ For $P_3$ we write
$$
P_3=\frac12\frac{\varpi}{|\da z|^2}(\da BR(z,\varpi)\cdot\dpa z+\da^2 z\cdot\dpa z),
$$
and we get
$$
|\dpt P_3|\leq C(\|\F(z)\|^2_{L^\infty}+\|z\|^2_{H^4}+\|\varpi\|^2_{H^3})^m(|\varpi_t|+|\da z_t|+|\da \varpi_t|+|\da^2 z_t|+|H(\da \varpi_t)|+|\la(\da z_t)|).
$$
It yields
$$\|\dpt P_3\|_{L^\infty}\leq C(\|\F(z)\|^2_{L^\infty}+\|z\|^2_{H^4}+\|\varpi\|^2_{H^3})^j(\|\varpi_t\|_{H^2}+ \|z_t\|_{H^3})$$
by the Sobolev embedding. The inequalities \eqref{nswt} and \eqref{nszt} take care of the rest.

In $\dpt P_4$ we have the term
$$\dpt \varphi=\frac{\varpi_t}{2|\da z|^2}-\varpi\frac{\da z\cdot\da z_t}{2|\da z|^3}-\dpt(|\da z|c)(\al,t),$$
but estimates \eqref{nswt} and \eqref{nszt} yield easily the
appropriate bounds for $\|\varphi_t\|_{L^\infty}$ and $\|\dpt
P_4\|_{L^\infty}$.

In a similar way we control $\|\dpt P_2\|_{L^\infty}$. Regarding
$\dpt P_1$ the most singular terms are given by
$$Q_1=\frac12H(\varpi_{tt}),\quad Q_2=-\frac1{2|\da z|^2} \la(z_{tt}\cdot\da z).$$
For $Q_2$ we  decompose further $Q_2 = R_1 + R_2$ where
$$
R_1=-\frac1{2|\da z|^2} H(z_{tt}\cdot\da^2 z),\quad R_2=-\frac1{2|\da z|^2} H(\da z_{tt}\cdot\da z).
$$
Then we take a time  derivative  in \eqref{eev} to estimate $R_1$ in
$L^\infty$,  and for $R_2$ we use the fact that $\da z_t\cdot\da z$
only depend on $t$ (see \eqref{AppA}). Next the identity $\da
z_{tt}\cdot\da z=\dpt(\da z_t\cdot\da z)-|\da z_t|^2$ allows us to
write
$$
R_2=\frac1{2|\da z|^2} H(|\da z_t|^2).
$$
>From estimates \eqref{nszt} we get  control of $R_2$ in $L^\infty$.

For $Q_1$ we have
$$
\|Q_1\|_{L^\infty}\leq C\|\varpi_{tt}\|_{C^\delta}.
$$
To continue we will need  estimates on $\|\varpi_{tt}\|_{C^\delta}$
for which we may use the identity \eqref{nose2}, and the inequality
$\|f\|_{C^\delta}\leq C(\|f\|_{L^2}+\|f\|_{\overline{C}^\delta})$
where
$$
\|f\|_{\overline{C}^\delta}=\sup_{\al\neq\beta}\frac{|f(\al)-f(\beta)|}{|\al-\beta|^\delta}.
$$
Then formula \eqref{nose2} gives
\begin{align}
\begin{split}\label{qt}
\varpi_{tt}+T(\varpi_{tt})&=\dpt I_1+\dpt I_2
-2\varphi_t\da\varphi-2\varphi\da\varphi_t+\dpt R+J_1+J_2,
\end{split}
\end{align}
where
$$
J_1=-\frac{1}{\pi}\int_{-\pi}^\pi \frac{(z_t(\al)\!-\!z_t(\al\!-\!\beta))^{\bot}\cdot\da
z(\al)}{|z(\al)-z(\al-\beta)|^2}\varpi_t(\al-\beta)d\beta,
$$
and
$$
J_2=\frac{2}{\pi}\int_{-\pi}^\pi \frac{(z(\al)\!-\!z(\al\!-\!\beta))^{\bot}\cdot\da
z(\al)}{|z(\al)-z(\al-\beta)|^4}(z(\al)\!-\!z(\al\!-\!\beta))\cdot(z_t(\al)\!-\!z_t(\al\!-\!\beta))\varpi_t(\al-\beta)d\beta.
$$
As before we use the invertibility of $(I+T)$ to get appropriate
estimates on $\|\varpi_{tt}\|_{L^2}$:
\begin{eqnarray}\label{nl2wtt}
\|\varpi_{tt} \|_{L^2}\leq
C\exp(C|||z|||^p)(\|\F(z)\|^2_{L^\infty}+\|z\|^2_{H^4}+\|\varpi\|^2_{H^3}+\|\varphi\|^2_{H^3})^j.
\end{eqnarray}
We shall show with some details how to get the most singular case
$\|\varpi_{tt}\|_{\overline{C}^\delta}$.

Formula \eqref{qt} yields
$$
\|\varpi_{tt}\|_{\overline{C}^\delta}\leq \|T(\varpi_{tt})\|_{\overline{C}^\delta}+\|\dpt I_1+\dpt I_2
-2\varphi_t\da\varphi-2\varphi\da\varphi_t+\dpt R+J_1+J_2\|_{\overline{C}^\delta},
$$
and therefore
$$
\|\varpi_{tt}\|_{\overline{C}^\delta}\leq \|T(\varpi_{tt})\|_{H^1}+\|\dpt I_1+\dpt I_2
-2\varphi_t\da\varphi-2\varphi\da\varphi_t+\dpt R+J_1+J_2\|_{H^1}.
$$
Then the inequality $\|T(\varpi_{tt})\|_{H^1}\leq
\|T\|_{L^2\rightarrow H^1}\|w_{tt}\|_{L^2}$, together with
\eqref{ol2h1} and \eqref{nl2wtt} yield the desired estimate. In
$\dpt I_1$ we find the term  $\la(z_{tt})$ therefore we need to
control $\|\la(z_{tt})\|_{H^1}=\|\da^2z_{tt}\|_{L^2}$, but formula
\eqref{eev} let us obtain that bound. In $\dpt I_2$ we have again
the extra cancelation given by
$$(z(\al)-z(\al-\beta))^\bot\cdot\da z(\al)=(z(\al)-z(\al-\beta)-\da
z(\al)\beta)^\bot\cdot\da z(\al),$$ which yields the appropriate
estimate. We have also to control $\|\da^2\varphi_t\|_{L^2}$, but
formula \eqref{phit} gives
$$
\da^2\varphi_t(\al,t)=\frac{\da^2\varpi_t(\al,t)}{2|\da
z(\al,t)|}-\frac{\da^2\varpi(\al,t)}{2|\da z(\al,t)|^{3}}\da
z(\al,t)\cdot\da z_t(\al,t)-\partial_t(\da(\da z\cdot\da
BR(z,\varpi))),
$$
showing that it can be estimated as before. Finally, the remainder
terms are less singular in derivatives, allowing us to finish the
proof.
%%%%%%%%%%%%%%%%%%%%%%%%%%%%%%%%%%%%%%%%%%%%%%%%%%%%%%%%%%%%%%%%%%%%%%%%%%%%%%%%%%%
%%%%%%%%%%%%%%%%%%%%%%%%%%%%%%%%%%%%%%%%%%%%%%%%%%%%%%%%%%%%%%%%%%%%%%%%%%%%%%%%%%%

\section{A priori energy estimates}

%%%%%%%%%%%%%%%%%%%%%%%%%%%%%%%%%%%%%%%%%%%%%%%%%%%%%%%%%%%%%%%%%%%%%%%%%%%%%%%%%%%
%%%%%%%%%%%%%%%%%%%%%%%%%%%%%%%%%%%%%%%%%%%%%%%%%%%%%%%%%%%%%%%%%%%%%%%%%%%%%%%%%%%

Let us consider for $k\geq 4$ the following definition of energy $E(t)$:
\begin{align}\label{E}
\begin{split}
E^2(t)&=\|z\|^2_{H^{k-1}}(t)+\int_{-\pi}^\pi\frac{\sigma(\al,t)}{\rho^2|\da z(\al,t)|^2}|\da^k z(\al,t)|^2 d\al\\
&\quad+\|\F(z)\|^2_{L^\infty}(t)+\|\varpi\|^2_{H^{k-2}}(t)+\|\varphi\|_{H^{k-\frac12}}^2(t),
\end{split}
\end{align}
so long as $\sigma(\al,t)>0$. In the next section we shall show a
proof of the following lemma.
\begin{lemma}
Let $z(\al,t)$ and $\varpi(\al,t)$ be a solution of (\ref{fibr}--\ref{fw}) in the case  $\rho_1=0$. Then, the following a priori
estimate holds:
\begin{align}
\begin{split}\label{ntni}
\frac{d}{dt}E^p(t)&\leq \frac{C}{m^q(t)}\exp(CE^p(t)),
\end{split}
\end{align}
for $m(t)=\D\min_{\al\in [-\pi,\pi]} \sigma(\al,t)=\sigma(\al_t,t)>0$, $k\geq 4$ and $C$, $q$ and $p$ some universal constants.
\end{lemma}

We shall present the details when $k=4$. Regarding  $\|\da^4
z\|^2_{L^2}$ let us remark that we have
$$
\|\da^4 z\|^2_{L^2}(t)=\int_\T
\frac{\sigma(\al,t)}{\sigma(\al,t)}|\da^4 z(\al,t)|^2d\al\leq
\frac{1}{m(t)} \int_\T \sigma(\al,t)|\da^4 z(\al,t)|^2d\al.
$$

%%%%%%%%%%%%%%%%%%%%%%%%%%%%%%%%%%%%%%%%%%%%%%%%%%%%%%%%%%%%%%%%%%%%%%%%%%%%%%%%%%%
%%%%%%%%%%%%%%%%%%%%%%%%%%%%%%%%%%%%%%%%%%%%%%%%%%%%%%%%%%%%%%%%%%%%%%%%%%%%%%%%%%%

\subsection{Energy estimates on the curve}

%%%%%%%%%%%%%%%%%%%%%%%%%%%%%%%%%%%%%%%%%%%%%%%%%%%%%%%%%%%%%%%%%%%%%%%%%%%%%%%%%%%
%%%%%%%%%%%%%%%%%%%%%%%%%%%%%%%%%%%%%%%%%%%%%%%%%%%%%%%%%%%%%%%%%%%%%%%%%%%%%%%%%%%

In this section we give the proof of the following  lemma when,
again, $k=4$. The case $k>4$ is left to the reader.

\begin{lemma}
Let $z(\al,t)$ and $\varpi(\al,t)$ be a solution of (\ref{fibr}--\ref{fw}) in the case  $\rho_1=0$.
 Then, the following a priori estimate holds:
\begin{align}
\begin{split}\label{eec}
\frac{d}{dt}\Big(\|z\|^2_{H^{k-1}}+\int_{-\pi}^\pi\frac{\sigma(\al)}{\rho^2|\da z(\al)|^2}|\da^k z(\al)|^2 d\al\Big)(t)&\leq S(t)+\frac{C}{m^q(t)}\exp(CE^p(t)),
\end{split}
\end{align}
for
\begin{equation}\label{fS}
S(t)=\int_{-\pi}^{\pi}\frac{2\sigma(\al)}{\rho^2}\frac{\da^k z(\al)\cdot \dpa z(\al)}{|\da
z(\al)|^3}\la(\da^{k-1}\varphi)(\al)d\al,
\end{equation}
and $k\geq 4$.
\end{lemma}
(We have denoted with $S$ a non integrable term which shall appear
in the equation of the evolution of $\varphi$ but with the opposite
sign.)\\

Proof: Using \eqref{nszt} one gets easily
\begin{align*}
\begin{split}
\frac{d}{dt}\|z\|^2_{H^3}&\leq C\int_{-\pi}^{\pi}(|z(\al)||z_t(\al)|+|\da^3z(\al)||\da^3z_t(\al)|)d\al\\
&\leq \frac{C}{m^q(t)}\exp(CE^p(t)).
\end{split}
\end{align*}
Then we have
\begin{align*}
\begin{split}
\frac{d}{dt} \int_{-\pi}^\pi\frac{\sigma(\al)}{\rho^2|\da z(\al)|^2}|\partial_{\al}^4 z(\al)|^2d\al
= &\int_{-\pi}^{\pi}\frac{1}{\rho^2}(\frac{\sigma_t(\al)}{|\da z(\al)|^2}-\frac{\sigma(\al) 2\da z(\al)\cdot\da z_t(\al)}{|\da z(\al)|^4})|\partial_{\al}^4 z(\al)|^2d\al\\
&+\int_{-\pi}^\pi\frac{2\sigma(\al)}{\rho^2|\da z(\al)|^2}\partial_{\al}^4 z(\al)\cdot\partial_{\al}^4 z_t(\al) d\al \\
=&I_1+I_2.
\end{split}
\end{align*}

The bound \eqref{nlinftyst} gives us
$$I_1\leq \frac{C}{m^q(t)}C\exp(C E^p(t)).$$

Next for $I_2$ we write
\begin{align*}
\begin{split}
I_2&=\int_{-\pi}^\pi\frac{2\sigma(\al)}{\rho^2|\da z|^2} \partial_{\al}^4 z(\al)\cdot
\partial_{\al}^4 BR(z,\varpi)(\al) d\al+
\int_{-\pi}^\pi\frac{2\sigma(\al)}{\rho^2|\da z|^2}\partial_{\al}^4 z(\al)\cdot \da^4(c\da z)(\al) d\al \\
&=J_1+J_2.
\end{split}
\end{align*}
The most singular terms in $J_1$ are given by $K_1$, $K_2$ and
$K_3$:

$$
K_1=\frac{1}{\pi}\int_{-\pi}^\pi\int_{-\pi}^\pi
\frac{\sigma(\al)}{\rho^2|\da z|^2}\partial_{\al}^4
z(\al)\cdot\frac{(\da^4z(\al)-\da^4z(\al-\beta))^\bot}{|z(\al)-z(\al-\beta)|^2}\varpi(\al-\beta)d\beta
d\alpha,$$

\begin{align*}
\begin{split}
K_2&=-\frac{1}{\pi}\int_{-\pi}^\pi\int_{-\pi}^\pi\frac{\sigma(\al)}{\rho^2|\da z|^2}\partial_{\al}^4
z(\al)\!\cdot\!\frac{(z(\al)\!-\!z(\al\!-\!\beta))^{\bot}}{|z(\al)\!-\!z(\al\!-\!\beta)|^4}C(\al,\beta)\varpi(\al\!-\!\beta)d\beta
d\alpha,
\end{split}
\end{align*}
and
$$
K_3=\frac{1}{\pi}\int_{-\pi}^\pi\int_{-\pi}^\pi
\frac{\sigma(\al)}{\rho^2|\da z|^2} \da^4z(\al)\cdot
\frac{(z(\al)-z(\al-\beta))^{\bot}}{|z(\al)-z(\al-\beta)|^2}\da^4\varpi(\al-\beta)
d\beta,
$$
where $C(\al,\beta)=(z(\al)-z(\al\!-\!\beta))\cdot(\da^4z(\al)\!-\!\da^4z(\al\!-\!\beta))$.

Then we write:
\begin{align*}
\begin{split}
K_1&=\frac{1}{\pi}\int_{-\pi}^\pi\int_{-\pi}^\pi
\frac{\sigma(\al)}{\rho^2|\da z|^2}
\partial_{\al}^4
z(\al)\cdot\frac{(\da^4z(\al)-\da^4z(\beta))^\bot}{|z(\al)-z(\beta)|^2}\varpi(\beta)d\beta
d\alpha\\
&=\frac{1}{\pi\rho^2|\da z|^2}\int_{-\pi}^\pi\int_{-\pi}^\pi\partial_{\al}^4
z(\al)\cdot\frac{(\da^4z(\al)-\da^4z(\beta))^\bot}{|z(\al)-z(\beta)|^2}\frac{\sigma(\al)\varpi(\beta)+\sigma(\beta)\varpi(\alpha)}{2}d\beta
d\alpha\\
&\quad+\frac{1}{\pi\rho^2|\da z|^2}\int_{-\pi}^\pi\int_{-\pi}^\pi\partial_{\al}^4
z(\al)\cdot\frac{(\da^4z(\al)-\da^4z(\beta))^\bot}{|z(\al)-z(\beta)|^2}\frac{\sigma(\al)\varpi(\beta)-\sigma(\beta)\varpi(\alpha)}{2}d\beta
d\alpha\\
&=L_1+L_2.
\end{split}
\end{align*}
That is we have performed a kind of integration by parts in $K_1$,
allowing us to show that $L_1$, its  most singular term,  vanishes:
\begin{align*}
\begin{split}
L_1&=-\frac{1}{\pi\rho^2|\da z|^2}\int_{-\pi}^\pi\int_{-\pi}^\pi\partial_{\al}^4
z(\beta)\cdot\frac{(\da^4z(\al)-\da^4z(\beta))^\bot}{|z(\al)-z(\beta)|^2}\frac{\sigma(\al)\varpi(\beta)+
\sigma(\beta)\varpi(\alpha)}{2}d\beta d\alpha\\
&=\frac{1}{2\pi\rho^2|\da z|^2}\int_{-\pi}^\pi\int_{-\pi}^\pi(\partial_{\al}^4
z(\al)\!-\!\partial_{\al}^4 z(\al))\!\cdot\!
\frac{(\da^4z(\al)\!-\!\da^4z(\beta))^\bot}{|z(\al)\!-\!z(\beta)|^2}\frac{\sigma(\al)\varpi(\beta)\!+\!
\sigma(\beta)\varpi(\alpha)}{2}d\beta d\alpha\\
&=0,
\end{split}
\end{align*}
whether for $L_2$ we have
\begin{align*}
\begin{split}
L_2&=-\frac{1}{\pi\rho^2|\da z|^2}\int_{-\pi}^\pi\int_{-\pi}^\pi\partial_{\al}^4
z(\al)\cdot\frac{(\da^4z(\beta))^\bot}{|z(\al)-z(\beta)|^2}
\frac{(\sigma(\al)-\sigma(\beta))\varpi(\beta)}{2}d\beta d\alpha\\
&\quad -\frac{1}{\pi\rho^2|\da z|^2}\int_{-\pi}^\pi\int_{-\pi}^\pi\partial_{\al}^4
z(\al)\cdot\frac{(\da^4z(\beta))^\bot}{|z(\al)-z(\beta)|^2}\frac{\sigma(\beta)(\varpi(\beta)-\varpi(\alpha))}{2}d\beta
d\alpha.
\end{split}
\end{align*}
In $L_2$ the kernels have degree $-1$ so long as the arc-chord
condition is satisfied, so they can be estimated by
$$L_2\leq C\|\F(z)\|^k_{L^\infty}\|z\|^k_{H^3}\|\varpi\|_{C^{1,\delta}}\|\sigma\|_{C^{1,\delta}}\|\da^4 z\|^2_{L^2}\leq \frac{C}{m^q(t)}\exp(C E^p(t)).$$

The term $C(\al,\beta)$ in $K_2$ can be written as follows:
\begin{align*}
C(\al,\beta)&=(z(\al)-z(\al\!-\!\beta)-\da z(\al)\beta)\cdot(\da^4z(\al)\!-\!\da^4z(\al\!-\!\beta))\\
  &\quad -\beta (\da z(\al)-\da z(\al\!-\!\beta))\cdot \da^4 z (\al\!-\!\beta)\\
  &\quad +\beta (\da z(\al)\cdot \da^4 z(\al)-\da z(\al\!-\!\beta)\cdot \da^4
  z(\al\!-\!\beta)),
\end{align*}
then using that
$$
\da z(\al)\cdot \da^4 z(\al)=-3\da^2 z(\al)\cdot \da^3 z(\al),
$$
we can split $K_2$ as a sum of kernels of degree $-1$ operating on
$\da^4z(\al)$, plus a kernel of degree $-2$ acting in three
derivatives $\da^3z(\al)$, allowing us to obtain again the estimate
$$K_2\leq \frac{C}{m^q(t)}\exp(C E^p(t)).$$
The term $K_3$ is a sum of a kernel of degree zero acting on four derivatives of $\varpi$

$$
L_3=\frac{1}{\pi}\int_{-\pi}^\pi\frac{\sigma(\al)}{\rho^2|\da z|^2} \da^4z(\al)\cdot\int_{-\pi}^\pi
[\frac{(z(\al)-z(\al-\beta))^{\bot}}{|z(\al)\!-\!z(\al\!-\!\beta)|^2}\!-\!\frac{\dpa z(\al)}{|\da z(\al)|^22\tan(\beta/2)}]\da^4\varpi(\al\!-\!\beta)
d\beta d\al,
$$
plus the following term:
$$
L_4=\int_{-\pi}^\pi\int_{-\pi}^\pi\frac{\sigma(\al)}{\rho^2}\frac{\da^4 z(\al)\cdot \da^{\bot}
z(\al)}{|\da z(\al)|^4}H(\da^4\varpi)(\al)
d\beta.
$$
We can integrate by parts on $L_3$ with respect to $\beta$ writing
$\da^4 \varpi(\al-\beta)=-\partial_{\beta}(\da^3\varpi(\al-\beta))$
and then pass  the derivative to the kernel of degree zero. This
calculation gives three derivatives in $\varpi$ and  kernels of
degree $-1$ which can be estimated as before.

 Next in  $L_4$ we write
\begin{align*}
\begin{split}
L_4&=\int_{-\pi}^\pi\frac{2\sigma(\al)}{\rho^2}\frac{\da^4 z(\al)\cdot \dpa z(\al)}{|\da z(\al)|^3} \la(\da^3(\frac{\varpi}{2|\da z|}))(\al) d\al\\
&=\int_{-\pi}^\pi\frac{2\sigma(\al)}{\rho^2}\frac{\da^4 z(\al)\cdot \dpa z(\al)}{|\da
z(\al)|^3} [\la(\da^3\varphi)(\al)-
\la(\da^2(\frac{\da z}{|\da z|}\cdot \da BR(z,\varpi))(\al)] d\al\\
&=S+M_0,
\end{split}
\end{align*}
for $S(t)$ given by \eqref{fS}. For $M_0$ we have
\begin{align*}
\begin{split}
\frac{\rho^2}{2}M_0&=\int_{-\pi}^\pi H\big(\sigma\frac{\da^4 z\cdot \dpa z}{|\da
z|^3}\big)(\al)\da^3(\frac{\da z}{|\da z|}
\cdot \da BR(z,\varpi))(\al) d\al=N_1+N_2+N_3,
\end{split}
\end{align*}
where
$$N_1=\int_{-\pi}^\pi H\big(\sigma\frac{\da^4 z\cdot \dpa z}{|\da z|^3}\big)(\al)\frac{\da^4 z}{|\da z|}
\cdot \da BR(z,\varpi)(\al) d\al,$$
$$N_2=\int_{-\pi}^\pi H\big(\sigma\frac{\da^4 z\cdot \dpa z}{|\da z|^3}\big)(\al)\frac{\da z}{|\da z|}
\cdot \da^4 BR(z,\varpi)(\al) d\al,$$ and $N_3$ is given by the rest
of the terms which can be controlled easily with the estimate that
we already have for the Birkhoff-Rott integral.

 Regarding  $N_1$ a
straightforward calculation gives
\begin{align*}
\begin{split}
N_1&\leq C\|\sigma\frac{\da^4 z\cdot \dpa z}{|\da
z|^3}\|_{L^2}\|\frac{\da^4 z}{|\da z|} \cdot \da
BR(z,\varpi)\|_{L^2}\\
&\leq C\|\sigma\|_{L^\infty}\|\F(z)\|^k_{L^\infty}\|\da
BR(z,\varpi)\|_{L^\infty}\|\da^4 z\|^2_{L^2}.
\end{split}
\end{align*}

Again, in $N_2$ we consider the most singular terms given by
$$
O_1=\int_{-\pi}^\pi H\big(\sigma\frac{\da^4 z\cdot \dpa z}{|\da
z|^3}\big)(\al)\frac{\da z(\al)}{|\da z(\al)|} \cdot
\frac{1}{2\pi}\int_{-\pi}^\pi \frac{(\da^4 z(\al)-\da^4
z(\al\!-\!\beta))^{\bot}}{|z(\al)-z(\al\!-\!\beta)|^2}\varpi(\al\!-\!\beta)
d\al,
$$
$$
O_2=-\int_{-\pi}^\pi H\big(\sigma\frac{\da^4 z\cdot \dpa z}{|\da
z|^3}\big)(\al)\frac{\da z(\al)}{|\da z(\al)|} \cdot
\frac{1}{4\pi}\int_{-\pi}^\pi
\frac{(z(\al)-z(\al\!-\!\beta))^{\bot}}{|z(\al)-z(\al\!-\!\beta)|^4}C(\al,\beta)\varpi(\al\!-\!\beta)
d\al,
$$
$$
O_3=\int_{-\pi}^\pi H\big(\sigma\frac{\da^4 z\cdot \dpa z}{|\da
z|^3}\big)(\al)\frac{\da z}{|\da z|} \cdot \da BR(z,\da^3
\varpi)(\al) d\al.
$$
Using the above decomposition for $C(\al,\beta)$ we can easily
estimate $O_2$. In $O_3$ we may write
$$\da z(\al)\cdot\da BR(z,\da^3 \varpi)(\al)=\frac12\da T(\da^3 \varpi)-\da^2z(\al)\cdot BR(z,\da^3 \varpi)(\al)$$
to obtain
$$\|\da z\cdot\da BR(z,\da^3 \varpi)\|_{L^2}\leq\|T(\da^3 \varpi)\|_{H^1}+\|\da^2z\|_{L^\infty}\|BR(z,\da^3 \varpi)\|_{L^2}$$
allowing us to control $O_3$.

Next we split $O_1$ into several kernels of degree one acting on
$(\da^4 z(\al))^{\bot}$, which can be estimated as before, plus the
term
$$
P_1=\frac{1}{2}\int_{-\pi}^\pi H\big(\sigma\frac{\da^4 z\cdot \dpa z}{|\da
z|^3}\big)(\al)\frac{\varpi(\al)\da z(\al)}{|\da z(\al)|^3} \cdot
\la((\da^4 z)^{\bot})(\al) d\al.$$

Then the following estimate for the commutator
$$\|\frac{\varpi\da z}{|\da z|^3} \cdot
\la((\da^4 z)^{\bot})-\la(\frac{\varpi\da z}{|\da
z|^3} \cdot (\da^4 z)^{\bot})\|_{L^2}\leq
\|\F(z)\|^3_{L^\infty}\|w\|_{H^2}\|z\|_{H^3}\|\da^4 z\|_{L^2},$$
yields
\begin{align*}
\begin{split}
P_1&\leq \|\F(z)\|^3_{L^\infty}\|w\|_{H^2}\|z\|_{H^3}\|\da^4
z\|_{L^2}-\frac12\int_{-\pi}^\pi \sigma\frac{\da^4 z\cdot \dpa z}{|\da
z|^3}\da \Big( \frac{\varpi  \da^4 z\cdot\dpa z}{|\da z|^3}
\Big)(\al) d\al
\end{split}
\end{align*}
using that
$$
\int_{-\pi}^\pi Hf(\al)\la g(\al)d\al=-\int_{-\pi}^\pi f(\al)\da
g(\al)d\al,
$$
and a straightforward integration by parts let us to control $P_1$.

So finally we have controlled  $J_1$ in the following manner:

$$
J_1\leq \frac{C}{m^q(t)}\exp(C E^p(t))+ S.
$$

To finish the proof  let us observe that the term $J_2$ can be
estimated integrating by parts,  using the identity $\da^4
z(\al,t)\cdot \da z(\al,t)=-3\da^3 z(\al,t)\cdot \da^2 z(\al,t)$ to
treat its most singular component. We have obtained
$$
\int_{\T}\frac{\sigma(\al)}{\rho^2|\da z|^2}\partial_{\al}^4 z(\al)\cdot\da
z(\al) \da^4c(\al) d\al =
3\int_{\T}\frac{1}{\rho^2|\da z|^2}\da(\sigma\partial_{\al}^3 z\cdot\da^2
z)(\al) \da^3c(\al) d\al
$$
and this yields the desired control. q.e.d.

%%%%%%%%%%%%%%%%%%%%%%%%%%%%%%%%%%%%%%%%%%%%%%%%%%%%%%%%%%%%%%%%%%%%%%%%%%%%%%%%%%%
%%%%%%%%%%%%%%%%%%%%%%%%%%%%%%%%%%%%%%%%%%%%%%%%%%%%%%%%%%%%%%%%%%%%%%%%%%%%%%%%%%%

\subsection{Energy estimates for the arc-chord condition}

%%%%%%%%%%%%%%%%%%%%%%%%%%%%%%%%%%%%%%%%%%%%%%%%%%%%%%%%%%%%%%%%%%%%%%%%%%%%%%%%%%%
%%%%%%%%%%%%%%%%%%%%%%%%%%%%%%%%%%%%%%%%%%%%%%%%%%%%%%%%%%%%%%%%%%%%%%%%%%%%%%%%%%%

In this section we analyze the evolution of the quantity
$\|\F(z)\|_{L^\infty}(t)$, which gives the local control of the
arc-chord condition.

\begin{lemma}\label{lemaarcchord} The following estimate holds
\begin{align}
\begin{split}\label{enlif}
\D\dt\|\F(z)\|^2_{L^\infty}(t)&\leq C(\|\F(z)\|^2_{L^\infty}(t)+\|z\|^2_{H^3}(t)+\|\varpi\|^2_{H^2}(t))^{j}.
\end{split}
\end{align}
\end{lemma}

Proof: First we compute the time derivative of the function $\F(z)$
as follows
\begin{align*}
\D\dt\F(z)(\al,\beta)(t)&=-\frac{|\beta|(z(\al,t)-z(\al-\beta,t))\cdot(z_t(\al,t)-z_t(\al-\beta,t))}{|z(\al,t)-z(\al-\beta,t)|^3},
\end{align*}
obtaining
$$
\D\dt\F(z)(\al,\beta)(t)\leq \frac{|\beta||z_t(\al,t)-z_t(\al-\beta,t)|}{|z(\al,t)-z(\al-\beta,t)|^2}\leq (\F(z)(\al,\beta)(t))^2\|\da z_t\|_{L^{\infty}}(t).
$$
Sobolev estimates and \eqref{nszt} yield
$$
\D\dt\F(z)(\al,\beta)(t)\leq C(\F(z)(\al,\beta)(t))^2
(\|\F(z)\|^2_{L^\infty}(t)+\|z\|^2_{H^{3}}(t)+\|\varpi\|^2_{H^{2}}(t))^j,
$$
and therefore
$$
\D\dt\F(z)(\al,\beta)(t)\leq C\F(z)(\al,\beta)(t)
\|\F(z)\|_{L^\infty}(t)(\|\F(z)\|^2_{L^\infty}(t)+\|z\|^2_{H^{3}}(t)+\|\varpi\|^2_{H^{2}}(t))^j,
$$
We shall denote $G(t)=C\|\F(z)\|_{L^\infty}(t)(\|\F(z)\|^2_{L^\infty}(t)+\|z\|^2_{H^{3}}(t)+\|\varpi\|^2_{H^{2}}(t))^j$, so that after an integration
in the time variable $t$ we get
\begin{align*}
\F(z)(t+h)&\leq \F(z)(t) \exp\,\big(
\!\int_t^{t+h}G(s)ds\big),
\end{align*}
and therefore
\begin{align*}
\|\F(z)\|_{L^\infty}(t+h)&\leq \|\F(z)\|_{L^\infty}(t) \exp\,\big(
\!\int_t^{t+h}G(s)ds\big),
\end{align*}
which yields
\begin{align*}
\D\dt\|\F(z)\|_{L^\infty}(t)&=\lim_{h\rightarrow
0^{+}}(\|\F(z)\|_{L^\infty}(t+h)-\|\F(z)\|_{L^\infty}(t))h^{-1}\\
&\leq \|\F(z)\|_{L^\infty}(t)\lim_{h\rightarrow 0^+}( \exp\,\big(
\!\int_t^{t+h}G(s)ds\big) -1)h^{-1}\leq \|\F(z)\|_{L^\infty}(t)G(t),
\end{align*}
allowing us to finish the proof of lemma \ref{lemaarcchord}. q.e.d.

%%%%%%%%%%%%%%%%%%%%%%%%%%%%%%%%%%%%%%%%%%%%%%%%%%%%%%%%%%%%%%%%%%%%%%%%%%%%%%%%%%%
%%%%%%%%%%%%%%%%%%%%%%%%%%%%%%%%%%%%%%%%%%%%%%%%%%%%%%%%%%%%%%%%%%%%%%%%%%%%%%%%%%%

\subsection{Energy estimates for $\varpi$ and $\varphi$}

%%%%%%%%%%%%%%%%%%%%%%%%%%%%%%%%%%%%%%%%%%%%%%%%%%%%%%%%%%%%%%%%%%%%%%%%%%%%%%%%%%%
%%%%%%%%%%%%%%%%%%%%%%%%%%%%%%%%%%%%%%%%%%%%%%%%%%%%%%%%%%%%%%%%%%%%%%%%%%%%%%%%%%%

In this section we complete the estimate \eqref{ntni} with the following result.

\begin{lemma}
Let $z(\al,t)$ and $\varpi(\al,t)$ be a solution of
(\ref{fibr}--\ref{fw}) in the case  $\rho_1=0$.  Then, the following
a priori estimate holds:
\begin{align}
\begin{split}\label{eec}
\frac{d}{dt}(\|\varpi\|^2_{H^{k-2}}+\|\varphi\|^2_{H^{k-\frac12}})(t)&\leq -S(t)+\frac{C}{m^q(t)}\exp(CE^p(t)).
\end{split}
\end{align}
for $k\geq 4$.
\end{lemma}

Proof: We shall present the details in the case $k=4$, leaving the
other cases to the reader.

Formula \eqref{nswt} shows easily that
$$
\frac{d}{dt}\|\varpi\|^2_{H^2}(t)\leq  (\exp(C|||z|||^p(t))(\|\F(z)\|^2_{L^\infty}(t)+\|z\|^2_{H^{4}}(t)+\|\varpi\|^2_{H^{3}}(t)+\|\varphi\|^2_{H^{3}}(t))^j)
$$
which together with \eqref{nsw} yields
$$
\frac{d}{dt}\|\varpi\|^2_{H^2}(t)\leq  \frac{1}{m^q(t)}C\exp(C
E^p(t)).
$$
Using the estimates obtained before one have
$$
\frac{d}{dt}\|\varphi\|^2_{L^2}(t)\leq  \frac{1}{m^q(t)}C\exp(C
E^p(t)).
$$
Next \eqref{eephi} yields
\begin{align}
\begin{split}
\frac{d}{dt}\|\la^{1/2}(\da^3 \varphi)\|^2_{L^2}(t)&=\int_\T \da^3\varphi(\al) \la(\da^3\varphi_t)(\al)  d\al=I_1+I_2+I_3+I_4,
\end{split}
\end{align}
where
$$
I_1=-\int_\T \frac{1}{2|\da
z|}\da^3\varphi(\al)\la(\da^4(\varphi^2))(\al)d\al,\qquad I_2=-\int_\T B(t)\da^3\varphi(\al)\la(\da^3\varphi)d\al,
$$
$$
I_3=-\int_\T\frac{1}{\rho^2|\da z|^3}\da^3\varphi(\al)\la(\da^2 (\sigma
\da^2z\cdot\dpa z))(\al)d\al,
$$
and
$$
I_4=-\int_\T\frac{1}{|\da z|^3}\da^3\varphi(\al)\la(\da^2(\da
BR(z,\varpi)\cdot \dpa z+\frac{\varpi}{2|\da z|^2}\da^2 z\cdot \dpa
z)^2)(\al)d\al.
$$

The most singular term in $I_1$ is given by
$$
J_1=-\int_\T \frac{1}{|\da
z|}\da^3\varphi(\al)\la(\varphi\da^4\varphi)(\al)d\al,
$$
and we have
$$
J_1=\int_\T \frac{1}{|\da z|}\la^{\frac12}(\da^3\varphi)(\al)
[\varphi(\al)\la^{\frac12}(\da^4\varphi)(\al)-\la^{\frac12}(\varphi\da^4\varphi)(\al)]d\al+
\int_\T \frac{\da \varphi(\al)}{2|\da
z|}|\la^{\frac12}(\da^3\varphi)(\al)|^2 d\al.
$$
The following estimate for the commutator $\|g\la^{\frac12}(\da
f)-\la^{\frac12}(g\da f)\|_{L^2}\leq \|g\|_{C^2}\|f\|_{H^{\frac12}}$
yields
$$
J_1\leq \|\F(z)\|_{L^\infty} \|\varphi\|^3_{H^{4-\frac12}},
$$
allowing us to get the estimate
$\D I_1\leq \frac{1}{m^q(t)}C\exp(C E^p(t)).$\\
The  boundedness of the term $B(t)$ gives us a similar control of
$I_2$
$$
I_2\leq \frac{1}{m^p(t)}C\exp(C E^p(t)).
$$
Next we write the term $I_4$ as follows:
$$
I_4=\int_\T\frac{1}{|\da z|^3}H(\da^3\varphi)(\al)\da^3(\da
BR(z,\varpi)\cdot \dpa z+\frac{\varpi}{2|\da z|^2}\da^2 z\cdot \dpa
z)^2(\al)d\al,
$$
where the most singular part is given by
$$
J_2=\int_\T\frac{2}{|\da
z|^3}H(\da^3\varphi)(\al)D(\al)\da^3D(\al)d\al,
$$
where
\begin{equation}\label{fD}
D(\al)=\da BR(z,\varpi)\cdot \dpa z+\frac{\varpi}{2|\da
z|^2}\da^2 z\cdot \dpa z.
\end{equation}
To analyze $\da^3(D)$, let us observe that the most singular terms
are given by
$$
E_1=\frac{1}{2\pi}\int_{-\pi}^\pi\frac{(\da^4 z(\al)-\da^4 z(\al-\beta)
)\cdot \da z(\al)}{|z(\al)-z(\al-\beta)|^2}\varpi(\al-\beta) d\al,
$$
$$
E_2=\frac{-1}{\pi}\int_{-\pi}^\pi\frac{(z(\al)-z(\al-\beta)
)\cdot \da z(\al)}{|z(\al)-z(\al-\beta)|^4}(z(\al)-z(\al-\beta)
)\cdot (\da^4 z(\al)-\da^4 z(\al-\beta))\varpi(\al-\beta) d\al,
$$
$$
E_3=BR(z,\da^4\varpi)\cdot \dpa z+\da^3(\frac{\varpi}{2|\da
z|^2}\da^2 z\cdot \dpa z).
$$
Since the terms $E_1$ and $E_2$ are singular only in the tangential
directions, we can again  use    the following identity
\begin{equation}\label{qg}
\da z(\al)\cdot \da^4 z(\al)=-3\da^2 z(\al)\cdot \da^3 z(\al),
\end{equation}
 to obtain the desired control.

 In $E_3$ the term $BR(z,\da^4\varpi)\cdot \dpa z$ can be written as the sum of
$\frac{1}{2}H(\da^4\varpi)$ plus kernels of degree zero in
$\da^4\varpi$, which are bounded in $L^2$. Therefore we can write it
as follows
$$
BR(z,\da^4\varpi)\cdot \dpa z=\frac{1}{2}H(\da^4\varpi)+\mbox{``bounded terms in $L^2$''}
$$

The identity
$$
\frac{1}{2}\da^4\varpi=|\da z|\da^4\varphi-\da^3(\da BR(z,\varpi)\cdot\da z)
$$
yields $$\da^3(\da BR(z,\varpi)\cdot\da
z)=H\big(\da^3(\frac{\varpi}{2|\da z|^2}(\da^2 z)^{\bot}\cdot \da
z)\big)+\mbox{``bounded terms in $L^2$''}.$$ That is
$$
\frac{1}{2}H(\da^4\varpi)=H(|\da z|\da^4\varphi)-
H^2\big(\da^3(\frac{\varpi}{2|\da z|^2}(\da^2 z)^{\bot}\cdot \da z)\big)+
\mbox{``bounded terms in $L^2$''}.
$$
and therefore
$$
\frac{1}{2}H(\da^4\varpi)=H(|\da z|\da^4\varphi)+\big(\da^3(\frac{\varpi}{2|\da z|^2}(\da^2 z)^{\bot}\cdot \da z)\big)+\mbox{``bounded terms in $L^2$''}.
$$
The above equality gives $E_3=|\da z|H(\da^4\varphi)+\mbox{``bounded terms in $L^2$''}.$

 Finally for $J_2$ we have
$$
J_2=\int_\T\frac{2}{|\da
z|^2}H(\da^3\varphi)(\al)H(\da^4\varphi)(\al)D(\al)d\al+\mbox{``bounded terms''},
$$
and an integration by parts gives us the desired estimate.

For $I_3$ it is important to arrange conveniently the derivatives
$$
I_3=-S+J_3+\mbox{``bounded terms''},
$$
where
\begin{equation}\label{J3}
J_3=\int_\T\frac{\da^2 z\cdot\dpa z}{|\da z|^3}H(\da^3\varphi)(\al)\da^3\sigma(\al)d\al.
\end{equation}
Then, because of its sign,  the  term involving the highest
derivative can be eliminated and we are left with the task of
estimating $J_3$. In order to do that we shall study the singular
term $\da^3\sigma(\al)$ using the splitting
\begin{align*}
\begin{split}
\da^3\sigma&=\da^3\big((\partial_t BR(z,\varpi)+\frac{\varphi}{|\da z|}\da BR(z,\varpi))\cdot \dpa z\big)\\
&\quad+\da^3\big(\frac12\frac{\varpi}{|\da z|^2}(\da
z_t+\frac{\varphi}{|\da z|}\da^2 z)\cdot
\dpa z\big)+\mathrm{g}\da^4 z_1\\
&=F_1+F_2+F_3.
\end{split}
\end{align*}
The term $F_3$ trivializes, whether  for $F_2$ we have
$$
\da^3\big(\frac12\frac{\varpi}{|\da z|^2}(\da BR(z,\varpi)\cdot
\dpa z+\frac{\varpi}{2|\da z|^2}\da^2 z\cdot
\dpa z\big)=\da^3\big(\frac12\frac{\varpi}{|\da z|^2}D \big)
$$
where $D$ is given by \eqref{fD} and the integral can be estimated
like $I_4$ or $J_2$. Finally we are left with $F_1$, and we shall
show that
\begin{equation}\label{off}
F_1=|\da z|H(\da^3\varphi_t)-c|\da z|H(\da^4\varphi)+\mbox{``bounded terms in $L^2$''}.
\end{equation}
Plugging the above decomposition in $J_3$ \eqref{J3} we can control
this term as before using the formula for $\da^3\varphi_t$
\eqref{eephi}.

Next we split $F_1=G_1+G_2$ where
$$
G_1=\da^3\big(\partial_t BR(z,\varpi)\cdot \dpa z\big),\quad G_2=
\da^3\big(\frac{\varphi}{|\da z|}\da BR(z,\varpi)\cdot \dpa z\big),
$$
and again we will consider the more singular terms. In $G_1$ we have
$$
O_1=\frac{1}{2\pi}\int_{-\pi}^\pi \frac{(\da^3 z_t(\al)-\da^3 z_t(\al-\beta))\cdot \da z(\al)}{|z(\al)-z(\al-\beta)|^2}\varpi(\al-\beta)d\beta,
$$
$$
O_2=-\frac{1}{\pi}\int_{-\pi}^\pi \frac{(z(\al)\!-\!z(\al\!-\!\beta))\cdot \da z(\al)}{|z(\al)\!
-\!z(\al\!-\!\beta)|^4}(z(\al)\!-\!z(\al\!-\!\beta))\!\cdot\!(\da^3z_t(\al)\!-\!\da^3z_t(\al\!-\!\beta)))\varpi(\al\!-\!\beta)d\beta,
$$
and
$$
O_3=\frac{1}{2\pi}\int_{-\pi}^\pi \frac{(z(\al)-z(\al-\beta))\cdot \da z(\al)}{|z(\al)-z(\al-\beta)|^2}\da^3\varpi_t(\al-\beta)d\beta.
$$
Let us write  $O_1=P_1+P_2$ where
$$
P_1=\frac{1}{2\pi}\int_{-\pi}^\pi \frac{\da^3 z_t(\al)\cdot \da z(\al)-\da^3
z_t(\al-\beta)\cdot \da z(\al-\beta)}{|z(\al)-z(\al-\beta)|^2}\varpi(\al-\beta)d\beta,
$$
and
$$
P_2=-\frac{1}{2\pi}\int_{-\pi}^\pi \frac{\da^3 z_t(\al-\beta)\cdot(\da z(\al)-\da z(\al-\beta)) }{|z(\al)-z(\al-\beta)|^2}\varpi(\al-\beta)d\beta.
$$
The term $P_2$ has a kernel of degree $-1$ in $\da^3 z_t$, giving us
a Hilbert
 integral of $\da^3 z_t$ which can be estimated using \eqref{nszt}. From its expression its follows that $P_1$
 can be written as the sum of   terms involving  kernels of degree $-1$ and the operator $\la$, that is:
  $$P_1=\frac{\varpi}{2|\da z|^2}\la(\da^3 z_t\cdot \da z)+\mbox{``bounded terms in $L^2$''}.$$
 Since $A'(t)=2\da z_t(\al,t)\cdot\da z(\al,t)$ we have
$$
\da^3 z_t\cdot \da z=-2\da^2z_t\cdot \da^2 z-\da z_t\cdot \da^3 z,
$$ which yields
$$P_1=\frac{\varpi}{2|\da z|^2}(-2\la(\da^2 z_t\cdot \da^2 z)-\la(\da z_t\cdot \da^3 z))
+\mbox{``bounded terms in $L^2$''}.$$ Then, as it was shown before,
the estimates for $z$ and $z_t$ give us the control of the term
$P_1$ in the $L^2$ norm.

Regarding $O_2$ we introduce into its integral expression the
following identity
\begin{align*}
(z(\al)\!-\!z(\al\!-\!\beta))\!\cdot\!(\da^3z_t(\al)\!-\!\da^3z_t(\al\!-\!\beta))&
=\beta\da z(\al)\!\cdot\!(\da^3z_t(\al)\!-\!\da^3z_t(\al\!-\!\beta))\\
+(z(\al)&\!-\!z(\al\!-\!\beta)-\da z(\al)\beta)\!\cdot\!(\da^3z_t(\al)\!-\!\da^3z_t(\al\!-\!\beta))
\end{align*} and then we just take the same steps that we followed with $O_1$.

Using the estimates \eqref{nswt} for $\varpi_t$ we get
$$O_3=\frac12H(\da^3 w_t)+\mbox{``bounded terms in $L^2$''},$$
and therefore
$$
G_1=\frac12H(\da^3 w_t)+\mbox{``bounded terms in $L^2$''}.
$$
The formula for $G_2$ gives us more singular terms, namely the
following ones
$$
O_4=\frac{\varphi}{2\pi|\da z|}\int_{-\pi}^\pi \frac{(\da^4 z(\al)
-\da^4 z(\al-\beta))\cdot \da z(\al)}{|z(\al)-z(\al-\beta)|^2}\varpi(\al-\beta)d\beta,
$$
$$
O_5=-\frac{\varphi}{\pi|\da z|}\int_{-\pi}^\pi \frac{(z(\al)\!-\!z(\al\!-\!\beta))
\cdot \da z(\al)}{|z(\al)\!-\!z(\al\!-\!\beta)|^4}(z(\al)\!-\!z(\al\!-\!\beta))\!
\cdot\!(\da^4z(\al)\!-\!\da^4z(\al\!-\!\beta)))\varpi(\al\!-\!\beta)d\beta,
$$
and
$$
O_6=\frac{\varphi}{2\pi|\da z|}\int_{-\pi}^\pi \frac{(z(\al)-z(\al-\beta))\cdot \da z(\al)}{|z(\al)-z(\al-\beta)|^2}\da^4\varpi(\al-\beta)d\beta.
$$
Using the identity \eqref{qg} we can estimate $O_4$ and $O_5$ as
before. Furthermore we have that $$O_6=\frac{\varphi}{2|\da
z|}H(\da^4\varpi)+\mbox{``bounded terms in $L^2$''},$$ and
$$G_2=\frac{\varphi}{2|\da z|}H(\da^4\varpi)+\mbox{``bounded terms in $L^2$''}.$$
Then we get
\begin{equation}\label{ufF1}
F_1=\frac12H(\da^3 w_t)+\frac{\varphi}{2|\da z|}H(\da^4\varpi)+\mbox{``bounded terms in $L^2$''}.
\end{equation}
We shall continue deducing \eqref{off}  from \eqref{ufF1} to
\eqref{off}, in order to do that let us write $$\frac12w_t=\dpt(|\da
z|)\frac{w}{2|\da z|}+|\da z|(\varphi_t+\dpt(|\da z| c))$$
$$
\frac12\da^3w_t=\dpt(|\da z|)\frac{\da^3w}{2|\da z|}+|\da z|\da^3\varphi_t-|\da z|\da^2\dpt(\frac{\da z}{|\da z|}\cdot \da BR(z,\varpi)).
$$
Since
\begin{align*}
|\da z|\da^2\dpt(\frac{\da z}{|\da z|}\cdot \da BR(z,\varpi))&=
\da^2(\da z\cdot \da\dpt BR(z,\varpi))+\da^2(\frac{\da z_t\cdot\dpa z}{|\da z|^2} \dpa z\cdot\da BR(z,\varpi))\\
&=\da^2(\da z\cdot \da\dpt BR(z,\varpi))+\mbox{``bounded terms in $L^2$''}.
\end{align*}
The last two identities allows us to consider
\begin{equation*}\label{dtwt}
\frac12\da^3w_t=|\da z|\da^3\varphi_t-\da^2(\da z\cdot \da\dpt BR(z,\varpi))+\mbox{``bounded terms in $L^2$''}
\end{equation*}
and therefore
\begin{equation*}\label{dtwt}
\frac12H(\da^3w_t)=|\da z|H(\da^3\varphi_t)-H(\da^2(\da z\cdot \da\dpt BR(z,\varpi)))+\mbox{``bounded terms in $L^2$''}
\end{equation*}
This formula indicates that to prove \eqref{off} it is enough to
obtain
\begin{align}
\begin{split}\label{off2}
c|\da z|H(\da^4\varphi)&=\frac{\varphi}{2|\da z|}H(\da^4\varpi)-G_3+\mbox{``bounded terms in $L^2$''},
\end{split}
\end{align}
where
\begin{equation}\label{G3}
G_3=H(\da^2(\da z\cdot \da\dpt BR(z,\varpi))).
\end{equation}
Again let us  consider the most singular terms in $\da^2(\da z\cdot \da\dpt BR(z,\varpi))$:
$$
O_7=\da(\da z\cdot BR(z,\da^2\varpi_t))$$
$$
O_8=\frac{-1}{\pi}\int_{-\pi}^\pi \frac{(z(\al)\!-\!z(\al\!-\!\beta))^{\bot}
\!\cdot\! \da z(\al)}{|z(\al)\!-\!z(\al\!-\!\beta)|^4}
(z(\al)\!-\!z(\al\!-\!\beta))\!\cdot\!(\da^3z_t(\al)\!-\!\da^3z_t(\al\!-\!\beta))\varpi(\al\!-\!\beta)d\beta,
$$
and
$$
O_9=\frac{1}{2\pi}\int_{-\pi}^\pi \frac{(\da^3z_t(\al)-\da^3z_t(\al-\beta))^{\bot}
\cdot \da z(\al)}{|z(\al)-z(\al-\beta)|^2}\varpi(\al-\beta)d\beta.
$$
The term $O_7=\frac12\da T(\da^2\varpi_t)$  is estimated in $L^2$
by using the operator $T$. In $O_8$ we substitute
$(z(\al)\!-\!z(\al\!-\!\beta))^{\bot}\cdot\da z(\al)$ by
$(z(\al)\!-\!z(\al\!-\!\beta)-\da z(\al)\beta)^{\bot}\cdot\da
z(\al)$ inside the integral
 and then we split the integral  in two terms ( one is multiplied by $\da^3 z_t(\al)$ and the other is an operator
  $R(\da^3 z_t)$ with kernel of degree $-1$) allowing us to
  integrate $O_8$.

  Regarding $O_9$ we have that
$$
O_9=\frac{-1}{2|\da z|^2}\la(\da^3z_t\cdot \dpa z\varpi)+\mbox{``bounded terms in $L^2$''},$$
and therefore the identity $H(\la(f))=-\da f$ yields for $G_3$ in \eqref{G3} the following configuration:
\begin{align*}
G_3&=\frac{1}{2|\da z|^2}\da(\da^3z_t\cdot \dpa z\varpi)+\mbox{``bounded terms in $L^2$''}\\
&=\frac{1}{2|\da z|^2}(\da(\da^3BR(z,\varpi)\cdot \dpa z\varpi)+\da(c\da^4 z\cdot \dpa z\varpi))+\mbox{``bounded terms in $L^2$''}\\
&=\frac{1}{2|\da z|^2}(\frac12H(\da^4\varpi)\varpi+c\varpi\da(\da^4 z\cdot \dpa z))+\mbox{``bounded terms in $L^2$''}.
\end{align*}
With this identity in \eqref{off2} we obtain
\begin{align*}
\frac{\varphi}{2|\da z|}H(\da^4\varpi)-G_3&=
-\frac{c}{2} H(\da^4\varpi)-\frac{c\,\varpi}{2|\da z|^2}\da(\da^4 z\cdot \dpa z)+\mbox{``bounded terms in $L^2$''}\\
&=-cH(|\da z|\da^4\varphi)-G_4+\mbox{``bounded terms in $L^2$''},
\end{align*}
for
\begin{equation*}\label{off3}
G_4=cH(|\da z|^2\da^4c)+\frac{c\,\varpi}{2|\da z|^2}\da(\da^4 z\cdot \dpa z).
\end{equation*}
Finally we only have to show that $G_4$ is a bounded term in $L^2$. But this follows because we have
$$
|\da z|^2\da^4c=-\da^3(\da z\cdot\da BR(z,\varpi))=\frac1{2|\da
z|^2}\la(\da^4 z\cdot\dpa z\varpi)+ \mbox{``bounded terms in
$L^2$''}.
$$

%%%%%%%%%%%%%%%%%%%%%%%%%%%%%%%%%%%%%%%%%%%%%%%%%%%%%%%%%%%%%%%%%%%%%%%%%%%%%%%%%%%%
%%%%%%%%%%%%%%%%%%%%%%%%%%%%%%%%%%%%%%%%%%%%%%%%%%%%%%%%%%%%%%%%%%%%%%%%%%%%%%%%%%%%

\section{The addition of the Rayleigh-Taylor condition to the energy}

%%%%%%%%%%%%%%%%%%%%%%%%%%%%%%%%%%%%%%%%%%%%%%%%%%%%%%%%%%%%%%%%%%%%%%%%%%%%%%%%%%%%
%%%%%%%%%%%%%%%%%%%%%%%%%%%%%%%%%%%%%%%%%%%%%%%%%%%%%%%%%%%%%%%%%%%%%%%%%%%%%%%%%%%%

Our final step is to use the a priori estimates to prove
local-existence (Theorem 1.1.). For that purpose we introduce a
regularized
 evolution equation which is well-posed  independently of the sign condition on $\sigma(\al,t)$ at $t=0$.  But for
$\sigma(\al,0)>0$, we shall find a time of existence  uniformly in
the regularization, allowing us to take the limit.

Let $z^\ep (\al,t)$ be a solution of the following system:
\begin{align*}
\begin{split}z^{\ep}_t(\al,t)&=BR(z^{\ep},\varpi^{\ep})(\al,t)+c^{\ep}(\al,t)\da z^{\ep}(\al,t),
\end{split}
\end{align*}
\begin{align*}
\begin{split}
\varpi^\ep_t&=-2\partial_t BR(z^\ep,\varpi^\ep)\cdot
\da z^\ep-\da((\varphi^\ep)^2)+2|\da z^\ep|B^\ep c^\ep-2\mathrm{g}\da z^\ep_2+\ep 2|\da z^\ep|\Delta\varphi^\ep,
\end{split}
\end{align*}
$z^\ep(\al,0)=z_0(\al)$ and $\varpi^\ep(\al,0)=\varpi_0(\al)$ for $\ep>0$, where

\begin{align*}
\begin{split}
c^{\ep}(\al)&=\frac{\al+\pi}{2\pi}\int_{-\pi}^\pi\frac{\da z^{\ep}(\al)}{|\da
z^{\ep}(\al)|^2}\cdot \da BR(z^{\ep},\varpi^{\ep})(\al) d\al-\int_{-\pi}^\al
\frac{\da z^{\ep}(\beta)}{|\da z^{\ep}(\beta)|^2}\cdot\partial_{\beta}
BR(z^{\ep},\varpi^{\ep})(\beta) d\beta,
\end{split}
\end{align*}
$$
\varphi^\ep=\D\frac{\varpi^\ep}{2|\da z^\ep|}-|\da z^\ep|c^\ep,\qquad B^\ep(t)=\frac{1}{2\pi}\int_{-\pi}^\pi \!\frac{\da z^\ep(\al,t)}{|\da
z^\ep(\al,t)|^2}\cdot \da BR(z^\ep,\varpi^\ep)(\al,t) d\al.
$$

Proceeding as in section 3 we find
\begin{align}
\begin{split}\label{fephitep}
\da\varphi^\ep_t&=-\frac{\da^2((\varphi^\ep)^2)}{2|\da
z^\ep|}-B^\ep(t)\da\varphi^\ep-\frac{\sigma^\ep}{\rho^2}\frac{\da^2z^\ep\cdot\dpa z^\ep}{|\da z^\ep|^3}-\partial_t(|\da z^\ep|B^\ep)\\
&\quad +\frac{1}{|\da z^\ep|^3}(\da BR(z^\ep,\varpi^\ep)\cdot \dpa
z^\ep+\frac{\varpi^\ep}{2|\da z^\ep|^2}\da^2 z^\ep\cdot \dpa z^\ep)^2+\ep\Delta\da\varphi^\ep,
\end{split}
\end{align}
where
\begin{align*}
\begin{split}
\frac{\sigma^\ep}{\rho^2}&=(\partial_t BR(z^\ep,\varpi^\ep)+\frac{\varphi^\ep}{|\da z^\ep|}\da BR(z^\ep,\varpi^\ep))\cdot \dpa z^\ep\\
&\quad+\frac12\frac{\varpi^\ep}{|\da z^\ep|^2}(\da
z^\ep_t+\frac{\varphi^\ep}{|\da z^\ep|}\da^2 z^\ep)\cdot
\dpa z^\ep+\mathrm{g}\da z^\ep_1.
\end{split}
\end{align*}
For this system there is local-existence for initial data
satisfying $\F(z_0)(\al,\beta)< \infty$ even if
$\sigma^\ep(\al,0)$ does not have the proper sign. In the
following we shall show briefly  how to obtain a solution of the
regularized system with $z^{\ep}, \varphi^\ep\in
C^1([0,\mathrm{T}^{\ep}],H^k)$ for $k\geq 4$. We shall prove the same a
priori estimates given in sections $6.1$, $6.2$ and $6.4$, but the
estimates corresponding to sections $6.3$ and $6.5$ are
respectively
\begin{align}
\begin{split}\label{ewtr}
\|\varpi^\ep_t\|_{H^k}&\leq C
\exp(C|||z^\ep|||^p)(\|\F(z^\ep)\|^2_{L^\infty}+\|z^\ep\|^2_{H^{k+2}}
+\|\varpi^\ep\|^2_{H^{k+1}}+\|\varphi^\ep\|^2_{H^{k+1}})^j\\
&\quad+\ep C\exp(C|||z^\ep|||^p)\|\Delta\da^k\varphi^\ep\|_{L^2} ,
\end{split}
\end{align}
\begin{align}
\begin{split}\label{esr}
\|\sigma^\ep \|_{H^k}&\leq
C\exp(C|||z^\ep|||^p)(\|\F(z^\ep)\|^2_{L^\infty}+\|z^\ep\|^2_{H^{k+2}}+\|\varpi^\ep\|^2_{H^{k+1}}+\|\varphi^\ep\|^2_{H^{k+1}})^j\\
&\quad +\ep C\exp(C|||z^\ep|||^p)\|\Delta\da^k\varphi^\ep\|_{L^2},
\end{split}
\end{align}
for $k\geq 2$.

Then following the same steps of section 6 we have
\begin{align*}
\frac{d}{dt}(\|z^\ep\|^2_{H^k}+&\|\F(z^\ep)\|^2_{L^\infty}+\|\varpi^\ep\|^2_{H^{k-2}}+\|\varphi^\ep\|^2_{H^k})^2(t)\\
&\leq
C(\ep)\exp ((\|z^\ep\|^2_{H^k}+\|\F(z^\ep)\|^2_{L^\infty}+
\|\varpi^\ep\|^2_{H^{k-2}}+\|\varphi^\ep\|^2_{H^k})^p(t))
\end{align*}
and where the only difference appears in the following new term
$$
I=-\int_{-\pi}^\pi\da^{k-1}\big(\frac{\sigma^\ep}{\rho^2}\frac{\da^2z^\ep\cdot\dpa
z^\ep}{|\da z^\ep|^3}\big) \da^k\varphi d\al\leq
C(\ep)\big\|\frac{\sigma^\ep}{\rho^2}\frac{\da^2z^\ep\cdot\dpa
z^\ep}{|\da
z^\ep|^3}\big\|^2_{H^{k-2}}+\ep\|\da^{k+1}\varphi^{\ep}\|^2_{L^2}
$$
which  is controlled by the Laplacian dissipation term introduced
in the regularization.

The next step is to integrate the system during a time $\mathrm{T}$
independent of $\ep$. We will show that for this system we have
\begin{align}
\begin{split}\label{ntniepsilon}
\frac{d}{dt}E^p(t)&\leq \frac{C}{(m^\ep)^q(t)}(\|\sigma_t^\ep\|_{L^\infty}+1)\exp(CE^p(t)),
\end{split}
\end{align}
where $E(t)$ is given by the analogous formula \eqref{E} for
the $\ep$-system,
$$m^\ep(t)=\D\min_{\al\in[-\pi,\pi]}\sigma^\ep(\al,t)=\sigma^\ep(\al_t,t)>0$$
and $C$, $p$ and $q$ universal constant independent of $\ep$.

  In
the following we shall select only the most singular terms,
showing for them the corresponding uniform estimates for $k=4$ and
leaving to the reader the remainder easier cases. Let us consider
the one corresponding to $I_3$ in section $7.3$, we have
$$
I^\ep_3=-\int_{\pi}^\pi\frac{1}{|\da z^\ep|^3}\da^3\varphi^\ep(\al)\la(\da^2(\sigma^\ep\da^2z^\ep\cdot\dpa z^\ep))d\al.
$$
We split this term as
$I^\ep_3=-S^\ep+J_2^\ep+J_3^\ep+\mbox{``bounded terms''}$ where
$S^\ep$ corresponds to $S$ in \eqref{fS} and
$$
J_2^\ep=\int_{\pi}^\pi\frac{C}{|\da z^\ep|^3}H(\da^3\varphi^\ep)(\al)\da^2(\sigma^\ep)\da(\da^2z^\ep\cdot\dpa z^\ep)d\al,
$$
$$
J_3^\ep=\int_{\pi}^\pi\frac{1}{|\da z^\ep|^3}H(\da^3\varphi^\ep)(\al)\da^3(\sigma^\ep)\da^2z^\ep\cdot\dpa z^\ep d\al.
$$
In $J_2^\ep$ we use \eqref{esr} to get
$$
J_2^\ep\leq \frac{C}{(m^\ep)^q(t)}\exp(CE^p(t))+\ep^2\|\da^4\varphi^\ep\|^2_{L^2}.
$$
The similarity with \eqref{J3} together with the use of the
corresponding version of  \eqref{off} allows us to get
$$
J_3^\ep=\int_{\pi}^\pi\frac{\da^2z^\ep\cdot\dpa z^\ep}{|\da z^\ep|^2}H(\da^3\varphi^\ep)H(\da^3\varphi^\ep_t) d\al
+M\ep^2\|\da^4\varphi^\ep\|^2_{L^2}+\mbox{``bounded terms''}
$$
that by formula \eqref{fephitep} becomes
$$
J_3^\ep=-\ep\int_{\pi}^\pi\frac{\da^2z^\ep\cdot\dpa z^\ep}{|\da z^\ep|^2}H(\da^3\varphi^\ep)H(\da^5\varphi^\ep) d\al
+M\ep^2\|\da^4\varphi^\ep\|^2_{L^2}+\mbox{``bounded terms''}.
$$
Then we can write it as follows
$$
J_3^\ep=-\ep\int_{\pi}^\pi\la^{\frac12}\Big(\frac{\da^2z^\ep\cdot\dpa
z^\ep}{|\da z^\ep|^2}
H(\da^3\varphi^\ep)\Big)\la^{\frac12}(\da^4\varphi^\ep) d\al
+M\ep^2\|\da^4\varphi^\ep\|^2_{L^2}+\mbox{``bounded terms''},
$$
and therefore
$$
J_3^\ep\leq M\ep^2\|\la^{\frac12}\da^4\varphi^\ep\|^2_{L^2}+\mbox{``bounded terms''}.
$$
Now the use of the Laplacian dissipative term introduced in the
evolution equation yields
\begin{align*}
\begin{split}
\frac{d}{dt}E^2(t)&\leq
\frac{C}{(m^\ep)^q(t)}(\|\sigma_t^\ep\|_{L^\infty}+1)\exp(CE^p(t))+
(M\ep^2-\ep)\|\la^{\frac12}\da^4\varphi^\ep\|^2_{L^2},
\end{split}
\end{align*}
where the constant $M$ is fixed. This finally shows \eqref{ntniepsilon} for $\ep$ small enough.

Our regularization damages the estimates for the term
$\|\sigma_t^\ep\|_{L^\infty}$ in \eqref{nlinftyst}.  But this
control is necessary only once in the argument and therefore
enough derivatives in the definition of energy gives the desired
control. Since we wish to keep the result for four derivatives, we
can go around the problem  just by regularizing the initial data.
At the end of the argument, when the local-existence theorem holds
for $\ep=0$, then the a priori energy estimate for $k=4$ allows us
to take the limit in the regularization of the initial data. With
this strategy and taking enough derivatives in the definition of
the energy, we find in \eqref{ntniepsilon} the following inequality
\begin{align}
\begin{split}\label{Epnd}
\frac{d}{dt}E^p(t)&\leq
\frac{C}{(m^\ep)^q(t)}\exp(CE^p(t)).
\end{split}
\end{align}

Now let us observe that if $z_0(\al)\in H^k$, $\varpi_0(\al)\in
H^{k-1}$ and $\varphi_0(\al)\in H^{k-\frac12}$, then we have the
solution in $[0,\mathrm{T}^\ep]$ of the regularized system. And if
initially $\sigma(\al,0)>0$, there is a time depending on $\ep$,
denoted by $\mathrm{T}^\ep$ again, in which $\sigma^\ep(\al,t)>0$. Now, for
$t\leq \mathrm{T}^\ep$ we have \eqref{Epnd}. Let us mention that at this
point of the proof we can not assume local-existence, because we
have the above estimate for $t\leq \mathrm{T}^\ep$, and if we let
$\ep\rightarrow 0$, it could be possible that $\mathrm{T}^\ep\rightarrow 0$
i.e. we cannot assume that if the initial data satisfy
$\sigma(\al,0)>0$,  there must be  a time $\mathrm{T}$, independent of
$\ep$, in which the following important quantity
$$m^\ep(t)=\D\min_{\al\in [-\pi,\pi]} \sigma^\ep(\al,t)=\sigma^\ep(\al_t,t)$$
is strictly grater that zero. In fact, everything in the evolution
problem depends upon the sign of $\sigma^\ep(\al,t)$  (the higher
order derivatives), since otherwise the problem is ill-posed
\cite{Ebin2}. In other words, at this stage of the proof  we do
not have local-existence when $\ep\rightarrow 0$, but the  following
argument will allow us to continue: First let us  introduce
the Rayleigh-Taylor condition in a new definition of energy as
follows:
$$
E_{RT}(t)=E^p(t)+\frac{1}{m^\ep(t)}.
$$
Sobolev inequalities shows that $\sigma^\ep(\al,t)\in
C^1([0,\mathrm{T}^\ep]\times[-\pi,\pi])$ and therefore $m^\ep(t)$ is a
Lipschitz function differentiable almost everywhere by
Rademacher's theorem.  With an analogous argument to the one used
in \cite{DY2} and \cite{DY3}, we can calculate the derivative of
$m^\ep(t)$, to obtain
$$
(m^\ep)'(t)=\sigma^\ep_t(\al_t,t)
$$
for almost every $t$. Then it follows that:
$$
\frac{d}{dt}\big(\frac{1}{m^\ep}\big)(t)=-\frac{\sigma^\ep_t(\al_t,t)}{(m^\ep)^2(t)}
$$
 almost everywhere. The control of the quantity
$\|\sigma_t^\ep\|_{L^\infty}$, independently of $\ep$,  by its
formula together with inequality \eqref{Epnd} yields
$$
\frac{d}{dt}E_{RT}(t)\leq C\exp(C E_{RT}(t)),
$$
and therefore
$$
E_{RT}(t)\leq -\frac{1}{C}\ln(\exp(-C E_{RT}(0)-C^2t),
$$
Now we are in position to  extend the time of existence $\mathrm{T}^\ep$ so
long as the above estimate works and obtain  a time $\mathrm{T}$
dependently only on the initial data (arc-chord and
Rayleigh-Taylor). Finally we can let $\ep$ tends to $0$ to
conclude the existence result.

%%%%%%%%%%%%%%%%%%%%%%%%%%%%%%%%%%%%%%%%%%%%%%%%%%%%%%%%%%%%%%%%%%%%%%%%%%%%%
%%%%%%%%%%%%%%%%%%%%%%%%%%%%%%%%%%%%%%%%%%%%%%%%%%%%%%%%%%%%%%%%%%%%%%%%%%%%%

\begin{quote}
\begin{tabular}{l}
\textbf{Antonio C\'ordoba} \\
{\small Departamento de Matem\'aticas}\\{\small Facultad de
Ciencias} \\ {\small Universidad Aut\'onoma de Madrid}
\\ {\small Crta. Colmenar Viejo km.~15,  28049 Madrid,
Spain} \\ {\small Email: antonio.cordoba@uam.es}
\end{tabular}
\end{quote}
\begin{quote}
\begin{tabular}{ll}
\textbf{Diego C\'ordoba} &  \textbf{Francisco Gancedo}\\
{\small Instituto de Ciencias Matem\'aticas} & {\small Department of Mathematics}\\
{\small Consejo Superior de Investigaciones Cient\'ificas} & {\small University of Chicago}\\
{\small Serrano 123, 28006 Madrid, Spain} & {\small 5734 University Avenue, Chicago, IL 60637}\\
{\small Email: dcg@imaff.cfmac.csic.es} & {\small Email: fgancedo@math.uchicago.edu}
\end{tabular}
\end{quote}


\begin{thebibliography}{99}


\bibitem{AM} D. Ambrose and N. Masmoudi. The zero surface tension limit of two-dimensional water waves.
\emph{Comm. Pure Appl. Math.} 58 1287-1315, 2005.

\bibitem{BMO} G. Baker, D. Meiron and S. Orszag. Generalized vortex methods for free-surface flow problems.
\emph{J. Fluid Mech.} 123 477-501, 1982.

\bibitem{BHL} T. Beale, T. Hou and J. Lowengrub. Growth rates for the linearized motion of fluid interfaces away from equilibrium.
 \emph{Comm. Pure Appl. Math.} 46, 1269-1301, 1993.
 
 \bibitem{ACG} A. Castro, D. C\'ordoba and F. Gancedo. A naive parametrization for the vortex-sheet problem. 
Preprint 2008, arXiv:0810.0731.

\bibitem{Chris} D. Christodoulou and H. Lindblad. On the motion of the free surface of a liquid.
\emph{Comm. Pure Appl. Math}. 53, no. 12, 1536--1602, 2000.

\bibitem{DY2} A. C\'ordoba, D. C\'ordoba and F. Gancedo. Interface evolution: the Hele-Shaw and Muskat problems.
Preprint 2008, arXiv:0806.2258.

\bibitem{DY3} A. C\'ordoba, D. C\'ordoba and F. Gancedo.  The Rayleigh-Taylor condition for the evolution of irrotational fluid interfaces.
Preprint, 2008.

\bibitem{DY} D. C\'ordoba and F. Gancedo. Contour dynamics of incompressible 3-D fluids in  a porous medium with different densities.
 \emph{Comm. Math. Phys.} 273, 2, 445-471 (2007).

\bibitem{Craig} W. Craig, An existence theory for water waves and the Boussinesq
and Kortewegde Vries scaling limits, \emph{Comm. Partial
Differential Equations}, 10, no. 8, 787-1003, 1985.

\bibitem{Shkoller} D. Coutand and S. Shkoller. Well-posedness of the free-surface incompressible Euler equations with or without surface tension.
\emph{J. Amer. Math. Soc.} 20, no. 3, 829--930, 2007.

\bibitem{Dahlberg} B.E.J. Dahlberg. On the Poisson integral for Lipschitz and
$C\sp{1}$-domains. \emph{Studia Math}. 66, no. 1, 13--24, 1979.

%\bibitem{Ebin1} D.G. Ebin. The equations of motion of a perfect fluid with free boundary are not well posed.
%\emph{Comm. Partial Differential Equations} 12, no. 10, 1175--1201, 1987.

\bibitem{Ebin2} D.G. Ebin. Ill-posedness of the Rayleigh-Taylor and Helmholtz problems for incompressible fluids.
\emph{Comm. Partial Differential Equations} 13, no. 10,
1265--1295, 1998.

\bibitem{Y} F. Gancedo. Existence for the $\alpha$-patch model and the QG sharp front in Sobolev spaces.
\emph{Adv. Math.}, Vol 217/6: 2569-2598, 2008.

\bibitem{Hou} T. Hou, J.S. Lowengrub and M.J. Shelley. Removing the Stiffness
from Interfacial Flows with Surface Tension. \emph{J. Comput. Phys.}, 114: 312-338, 1994.

\bibitem{Lannes} D. Lannes, Well-posedness of the water-waves equations,
\emph{J. Amer. Math. Soc.}, 18,  605-654, 2005.

\bibitem{Lindblad} H. Lindblad. Well-posedness for the motion of an incompressible liquid with free surface boundary.
 \emph{Ann. of Math.} (2) 162, no. 1, 109--194, 2005.

\bibitem{Nalinov} V.I. Nalinov. The Cauchy-Poisson Problem (in Russian),
 \emph{Dynamika Splosh. Sredy}, 18,1040-210, 1974.

\bibitem{Ray}Lord Rayleigh (J.W. Strutt), On the instability of jets.
\emph{Proc. Lond. Math. Soc}. 10, 4-13, 1879.

\bibitem{Shatah} J. Shatah and C. Zeng. Geometry and a priori estimates for free boundary problems of the Euler equation.
\emph{Comm. Pure Appl. Math.} 61, no. 5, 698--744, 2008.

\bibitem{St3} E.~Stein. \newblock Harmonic Analysis. \newblock \emph{
Princeton University Press.} Princeton, NJ, 1993.

\bibitem{Taylor} G. Taylor. The instability of liquid surfaces when accelerated in a direction perpendicular
to their planes. I. \emph{Proc. Roy. Soc. London. Ser. A.} 201,
192-196, 1950.

\bibitem{Wu} S. Wu. \newblock Well-posedness in Sobolev spaces of the full water wave problem in 2-D.
\emph{Invent. math.} 130, 39-72, 1997.

\bibitem{Wu2} S. Wu. \newblock Well-posedness in Sobolev spaces of the full water wave problem in 3-D.
\emph{J. Amer. Math. Soc.} 12, 445-495, 1999.

\bibitem{Yoshira} H. Yosihara, Gravity Waves on the Free Surface of an
Incompressible Perfect Fluid, \emph{Publ. RIMS Kyoto Univ.}, 18,
49-96, 1982.

\bibitem{ZZ} P. Zhang and Z. Zhang. On the free boundary problem of three-dimensional incompressible euler equations.
\emph{Comm. Pure and Appl. Math.}, 61: 877-940, 2008.

\end{thebibliography}
\end{document}